# Structure and Stability of Certain Chemical Networks and Applications to the Kinetic Proofreading Model of T-Cell Receptor Signal Transduction

Eduardo D. Sontag*

October 26, 2018

**Abstract**

This paper deals with the theory of structure, stability, robustness, and stabilization for an appealing class of nonlinear systems which arises in the analysis of chemical networks. The results given here extend, but are also heavily based upon, certain previous work by Feinberg, Horn, and Jackson, of which a self-contained and streamlined exposition is included. The theoretical conclusions are illustrated through an application to the kinetic proofreading model proposed by McKeithan for T-cell receptor signal transduction.

## 1 Introduction

This work was originally motivated by the study of the following system of first-order ordinary differential equations:

$$\dot{C}_0 = k_1 \left(T^* - \sum_{i=0}^{N} C_i\right)\left(M^* - \sum_{i=0}^{N} C_i\right) - (k_{-1,0} + k_{p,0})\, C_0$$

$$\vdots$$

$$\dot{C}_i = k_{p,i-1} C_{i-1} - (k_{-1,i} + k_{p,i})\, C_i$$

$$\vdots$$

$$\dot{C}_N = k_{p,N-1} C_{N-1} - k_{-1,N} C_N$$

where the subscripted $k$'s, as well as $M^*$ and $T^*$, are arbitrary positive constants; the $C_i$'s are nonnegative functions of time $t$, and dots indicate derivatives with respect to $t$. These equations arise in immunology, and describe a possible mechanism, due to McKeithan, cf. [19], that may explain the selectivity of T-cell interactions (some more details are provided below). McKeithan analyzed the equilibria of these equations, which represent steady-state regimes, mostly under the simplifying assumptions that $k_{p,i} \equiv k_p$ and $k_{-1,i} \equiv k_{-1}$ for some fixed $k_p$ and $k_{-1}$. During Carla Wofsy's series of talks [24] on the topic, a number of questions arose: what can be said about the structure of equilibria? are the equilibria stable? Stability is a natural mathematical requirement, with clear biological significance.

Under the simplifying assumptions mentioned above, it is not difficult to see that there is only one equilibrium, whose coordinates depend, of course, on the given constants. We will see that in fact even in the general case of the constants being distinct, there is a *unique equilibrium*. Far more interestingly, however, we will prove that the equilibrium is *globally asymptotically stable:* every solution, for any choice of initial (nonnegative) values $C_i(0)$'s, converges to the unique equilibrium. This conclusion rules out, in particular, periodic orbits and, of course, chaotic behaviors, and shows the "determinism" of the process described by McKeithan. Moreover, we will also establish the *robustness* of stability with respect to a quantifiable class of perturbations in the dynamics.

---

*Eduardo Sontag is with the Department of Mathematics, Rutgers University, New Brunswick, NJ 08903, USA. E-mail: sontag@hilbert.rutgers.edu, URL: http://www.math.rutgers.edu/ sontag This work was supported in part by US Air Force Grant F49620-98-1-0242.

Going further, and approaching now the equations from the point of view of a control theorist, one may pose questions of a very different nature as well, questions of *design*. The equilibria depend, in particular, on the value of the constants $T^*$ and $M^*$, which represent total concentrations (of intermediate complexes plus T-cell receptors, and intermediates plus peptide-major histocompatibility complex, respectively). If these concentrations are allowed to vary, and are seen as additional state variables, multiple attracting equilibria exist. So one may ask: if we are allowed to manipulate some of the variables, can we change the equilibria at will (and preserve stability)? How many variables need to be manipulated? Such questions might eventually impact approaches to therapy and rational drug design. We will provide a global feedback stabilization result as an answer.

As it turns out, the system of interest can be viewed as a very special instance of a very large class of nonlinear systems, for which the above-mentioned results can be established in general. *In this paper, we describe a theory of structure, stability, and stabilization for that general class of systems.* Thus, this paper can be read totally independently of the above-mentioned motivating example from immunology. By providing a general theory, one expects that other applications of the theorems given here will be possible.

The class of systems which we consider is, basically, that of "deficiency zero chemical reaction networks with mass-action kinetics (and one linkage class)" in the language of the beautiful and powerful theory developed by Feinberg, Horn, and Jackson, cf. [9, 10, 11, 12, 16]. As a matter of fact, our stability (but not the robustness nor the feedback) results are basically contained in that previous work. The existence and uniqueness of equilibria, and local asymptotic stability, are already proved in the papers by Feinberg, Horn, and Jackson, and we claim no originality whatsoever in that regard. And although the *global* stability results may not be readily apparent from a casual reading of that literature, our proofs of them consist basically of a careful repackaging of the discussion found in Feinberg's paper [9]. (One should point out that an alternative approach to global stability, which would apply to a somewhat smaller class of systems, could use ideas from [20], which, in turn, was based on [23]. Also somewhat related, but far more restrictive, are the results described in the book [14], which apply to systems which satisfy a "conservation of energy" constraint.)

The reliance upon the Feinberg/Horn/Jackson theory notwithstanding, we provide here a totally self-contained and streamlined exposition even of those results that are known for these systems (deep theories for other classes of systems have also been developed by these authors). We do this in order to present things in a terminology and formalism more standard in control and dynamical systems, and also because it is not obvious how to put together the necessary pieces from the various sources in a way that allows us to refer to them efficiently. In any case, it is our hope that this exposition will serve to make a wider audience in the dynamical systems and control theory communities aware of their work.

In addition to the stability results, we also prove a robustness and a feedback stabilization result, and we show the existence of a global change of coordinates, which brings the systems in question into a canonical form which exhibits a particularly elegant structure. We develop a formalism, and present explicit estimates for stability margins, with a view towards further theoretical developments. Indeed, in a recent follow-up joint paper [6] with Madalena Chaves, we have been able to derive, using these techniques, various input-to-state stability results and, in particular, a design of globally convergent Luenberger-like observers for systems of the type studied here. As new results, we also include partial generalizations to non-mass action kinetics (cf. Remark 5.6 and the stability section).

There is an extensive literature regarding applications of control theory to chemical engineering, and, in particular, to bioreactors. Among many references, one may mention the textbook [2] and the more recent survey articles [4] and [1], and the references given there, as well as the well-known theory of compartmental systems (see e.g. [18]). More specifically, results based on Feinberg-Horn-Jackson theory have appeared in the control literature, see [8, 5]. Other reactor biocontrol work deals with reachability and controllability issues, see for instance [3]. A preliminary version of this paper was posted electronically in [22].

The organization of this paper is as follows. In the rest of this Introduction, we motivate the formalism and results to follow by reviewing some basic facts concerning chemical reaction networks and working through an example. This section can be skipped with no loss of continuity by readers interested in the mathematical developments. The paper starts in Section 2, where we introduce the class of dynamical systems being studied, and the main theorems are also stated. In Section 3, we specialize to McKeithan's system, interpreting the results in that special case. Section 4 discusses some basic coordinatization facts, as well as some useful alternative system descriptions. Section 5 deals with the main proofs regarding interior equilibria, Section 6 with boundary equilibria, Section 7 with technical invariance results, and, finally, Section 8 with the proofs of the stability theorems.

## Acknowledgments

The author wishes to express great thanks to Carla Wofsy for a fascinating series of lectures, to Leah Keshet for organizing a superb workshop, and most especially to Marty Feinberg for making available reprints of his work and for very enlightening e-mail discussions.

## 1.1 Chemical Networks

Let us motivate the formalism and results to follow by reviewing some basic facts concerning chemical reaction networks. We will restrict attention to "mass action kinetics"; however, it is important to remark that variations such as Michaelis-Menten reactions, obtained through singular perturbation analysis when starting from mass-action models, are routinely used to model many enzymatic reactions in biology.

In chemical network studies, one analyzes systems of differential equations which describe the time-evolution of the concentrations $x_1(t), \ldots, x_n(t)$ of $n$ given "chemical species" $P_1, \ldots, P_n$; the $P_i$'s might denote anything from small molecules to large complexes. The equations are derived from a consideration of the reactions that are known to occur among the substances $P_i$, perhaps helped by other molecules which are not explicitly considered in the equations (such as catalysts or energy sources).

As an extremely simple illustration, suppose that each molecule of a certain species $P_1$ can react with one molecule of $P_2$ to produce a molecule of $P_3$, and that, conversely, each molecule of $P_3$ may dissociate (through a different process, typically with different time constants), into $P_1$ and $P_2$. This is indicated graphically by

$$P_1 + P_2 \to P_3 \quad \text{and} \quad P_3 \to P_1 + P_2,$$

or just by

$$P_1 + P_2 \rightleftharpoons P_3.$$

An example of this reaction is provided by the synthesis of ethyl tert-butyl ether ($P_3 = C_6H_{14}O$) from isobutene ($P_1 = C_4H_8$) and ethanol ($P_2 = C_2H_6O$). For studies of the ethyl tert-butyl ether reaction in the control theory literature, see [8], and [13]. Assuming that the reactor is well-mixed, particle-collision theories or quantum-mechanical potential energy methods are often used to justify the statement that the probability of such a reaction occurring, in a small time interval around time $t$, is proportional to the product $x_1(t)x_2(t)$ of the concentrations at time $t$, and to the length of the interval, that is to say, to the probability of two molecules in this bimolecular reaction colliding by virtue of being "in the same place at the same time". Since we gain a single molecule of $P_3$ for each such reaction, we arrive at a formula for the rate of increase of the concentration of $P_3$ due to the first reaction $P_1 + P_2 \to P_3$:

$$\dot{x}_3 = k_1 x_1 x_2. \tag{1}$$

In addition, as one molecule of $P_1$ and $P_2$ each are eliminated at the same rate, we also have the following two differential equations:

$$\dot{x}_1 = -k_1 x_1 x_2, \quad \dot{x}_2 = -k_1 x_1 x_2. \tag{2}$$

Here, $k_1$ is a suitable constant of proportionality, the "reaction rate constant", which is often taken to be (Arrhenius law) proportional to the Boltzmann factor $e^{-E/RT}$ and $T$ is the temperature. One also writes graphically:

$$P_1 + P_2 \xrightarrow{k_1} P_3.$$

We also assumed for our example that there is a dissociation reaction, that is,

$$P_3 \xrightarrow{k_2} P_1 + P_2$$

where $k_2$ is another rate constant, and this reaction is also modeled by rate equations: $P_3$ decays at a rate proportional to its concentration, $\dot{x}_3 = -k_2 x_3$, and each of $x_1$ and $x_2$ grow at this rate. Incorporating these into the previous equations gives the final set of differential equations:

$$\begin{align}
\dot{x}_1 &= -k_1 x_1 x_2 + k_2 x_3 \\
\dot{x}_2 &= -k_1 x_1 x_2 + k_2 x_3 \\
\dot{x}_3 &= k_1 x_1 x_2 - k_2 x_3,
\end{align} \quad (3)$$

which describe the evolution of all the concentrations $x_i(t)$.

A very convenient and systematic formalism to describe the complete system of equations is as follows. The entire reaction is represented by a graph, whose nodes are the "complexes" which appear in the reactions, such as $P_1 + P_2$ and $P_3$ in the example given above, and whose edges are labeled by the reaction rate constants. So, in the example shown above, there is an edge labeled $k_1$ (where $k_1$ is an actual positive number) starting at the node in the graph corresponding to $P_1 + P_2$, and pointing to the node corresponding to $P_3$, and there is likewise an edge labeled $k_2$ from the second node to the first. We associate to this graph its $2 \times 2$ incidence (connectivity) matrix matrix $A = \{a_{ij}\}$, listing all the edge labels (for instance, $a_{21} = k_1$, to indicate a reaction with rate constant $k_1$, from the first node to the second node). More generally, the size of $A$ is $m \times m$, if there are a total of $m$ complexes. One of the hypotheses to be made is that the matrix $A$ is irreducible, meaning that the graph is strongly connected, that is, there is some path, typically through several stages of intermediate reactions, linking any two given complexes (in fact, results hold under a weaker, block-irreducibility, property). Next, one introduces a set of column $n$-vectors $b_1, \ldots, b_m$, one for each complex (in our simple example, $m = 2$). This is done by specifying the contributions from each type of molecule. For example, $P_1 + P_2$ gives rise to the vector $b_1 = (1, 1, 0)'$, and $P_3$ to the vector $b_2 = (0, 0, 1)'$.

The "mass action" dynamics are then summarized by the system

$$\dot{x} = \sum_{i=1}^{m} \sum_{j=1}^{m} a_{ij}\, x_1^{b_{1j}} x_2^{b_{2j}} \ldots x_n^{b_{nj}} (b_i - b_j) \quad (4)$$

(where each $b_\ell \in \mathbb{R}^n$ has entries $b_{1\ell}, \ldots, b_{n\ell}$), which the reader may easily verify reduces to (3) in the above example, for which

$$A = \begin{pmatrix} 0 & k_2 \\ k_1 & 0 \end{pmatrix}, \quad b_1 = \begin{pmatrix} 1 \\ 1 \\ 0 \end{pmatrix}, \quad b_2 = \begin{pmatrix} 0 \\ 0 \\ 1 \end{pmatrix}.$$

Note that, for instance, the term corresponding to $i = 2$ and $j = 1$ in the summation in (4) gives us

$$k_1 x_1^0 x_2^1 x_3^1 (b_2 - b_1) = \begin{pmatrix} -k_1 x_1 x_2 \\ -k_1 x_1 x_2 \\ k_1 x_1 x_2 \end{pmatrix}$$

which provides precisely the contributions to $\dot{x}_1$, $\dot{x}_2$, and $\dot{x}_3$ represented by (1) and (2). (Higher-order polynomial equations may result as well. For example, suppose that each molecule of species $P_1$ can react with four molecules of $P_2$ to produce two molecules of $P_3$, that is, $P_1 + 4P_2 \to 2P_3$. This would give us equations $\dot{x}_1 = -k x_1 x_2^4$, $\dot{x}_2 = -4k x_1 x_2^4$, and $\dot{x}_3 = 2k x_1 x_2^4$, and the vectors $b_1$ and $b_2$ would become $b_1 = (1, 4, 0)'$, $b_2 = (0, 0, 2)'$.)

A fundamental role is played by the linear subspace $\mathcal{D}$ of $\mathbb{R}^n$ which is spanned by all the differences $b_i - b_j$. This is the *stoichiometric subspace* associated to the reaction, and each intersection between a parallel translate of $\mathcal{D}$ and the positive orthant is called a *class*, or more properly a stoichiometric compatibility class ("stoicheion" = (Gk.) element). The significance of classes is that, since $\dot{x} \in \mathcal{D}$, trajectories remain in classes, that is, classes are positive-time invariant manifolds for the dynamical system. (The positive orthant is itself forward invariant, as is easily shown.) In the example discussed above, $\mathcal{D}$ is the line spanned by $b_1 - b_2 = (1, 1, -1)'$. This line can also be described as the set of solutions of $x_1 + x_3 = x_2 + x_3 = 0$, so each class is given by the positive points in $\{x_1 + x_3 = c_1, x_2 + x_3 = c_2\}$, for different constants $c_1$ and $c_2$. Of course, it is clear from the equations (3) that $d(x_1 + x_3)/dt = d(x_2 + x_3)/dt \equiv 0$ along solutions.

One of the basic facts about the systems studied here is that there is a unique equilibrium in each class, and, under mild conditions, this equilibrium is globally asymptotically stable with respect to positive initial conditions. Continuing with the above example, the set of possible equilibria consists of the points in the hyperbolic paraboloid $k_2 x_3 - k_1 x_1 x_2 = 0$. Let's now take for simplicity $k_1 = k_2 = 1$ and analyze the stability of the equilibrium $x^0 = (1, 1, 1)'$. As the ray $(1, 1, 1)' + \lambda(1, 1, -1)$ is forward invariant, where $\lambda$ ranges over the interval $(-1, 1)$ (the intersection of the corresponding line and the positive orthant), we may parameterize motions by $\lambda$. One obtains the scalar differential equation $\dot{\lambda} = -3\lambda - \lambda^2$, which has, as claimed, the point $\lambda = 0$ as an asymptotically stable state with a basin of attraction which includes all of $(-1, 1)$.

If we start at a point which is not in the above line, the equilibrium approached will not be $x^0$. If this equilibrium $x^0$ is desirable, we may want to design a feedback law to drive the solutions from every other (positive orthant) initial state into $x^0$. Suppose that we can control the inflows and outflows of, let us say, $P_1$ and $P_3$. This situation is represented by the following control system:

$$\dot{x}_1 = -x_1 x_2 + x_3 + u_1, \quad \dot{x}_2 = -x_1 x_2 + x_3, \quad \dot{x}_3 = x_1 x_2 - x_3 + u_2.$$

Stabilization of $x^0 = (1, 1, 1)'$ can be achieved for instance by a simple feedback linearization, taking $u_1 := x_1 x_2 - x_3 - x_1 + 1$ and $u_2 := -x_1 x_2 + 1$. Then, $x_1(t) \to 1$ and $x_3(t) \to 1$ as $t \to \infty$, and thus, by a cascade argument, we see that also $x_2(t) \to 1$ as $t \to \infty$, as wanted. We will provide a general result on stabilization (but using linear feedback), as well as results on robustness of stability and on global decompositions of dynamics, for systems (4). These results will be illustrated again for the above simple example, now viewed as a particular case of McKeithan's system, in Section 3.

## 2 Definitions and Statements of Main Results

Some standard notations to be used are:

- $\mathbb{R}_{\geq 0}$ (resp., $\mathbb{R}_+$) = nonnegative (resp., positive) real numbers
- $\mathbb{R}_+^n$ (resp., $\mathbb{R}_+^{m \times m}$) = $n$-column vectors (resp., $m \times m$ matrices) with entries on $\mathbb{R}_+$; similarly for $\mathbb{R}_{\geq 0}$
- $\mathbb{R}_0^n$ = boundary of $\mathbb{R}_{\geq 0}^n$, set of vectors $x \in \mathbb{R}_{\geq 0}^n$ such that $x_i = 0$ for at least one $i \in \{1, \ldots, n\}$
- $x'$ = transpose of vector or matrix $x$
- $|x|$ = Euclidean norm of vector in $\mathbb{R}^n$
- $\langle x, z \rangle = x'z$, inner product of two vectors
- $\mathcal{D}^\perp = \{x \,|\, \langle x, z \rangle = 0 \,\forall z \in \mathcal{D}\}$.

Although we develop the theory for a somewhat wider class of systems, we wish to emphasize that *all the results to be given are valid, in particular, for the following general class of systems evolving on* $\mathbb{R}_{\geq 0}^n$:

$$\dot{x} \;=\; \sum_{i=1}^{m}\sum_{j=1}^{m} a_{ij}\, x_1^{b_{1j}} x_2^{b_{2j}} \ldots x_n^{b_{nj}}\, (b_i - b_j)$$

(and, as appropriate, for those systems obtained by adding control inputs, as discussed below). Each column vector $b_\ell \in \mathbb{R}^n$ has entries $b_{1\ell}, \ldots, b_{n\ell}$, which are nonnegative integers, and the $a_{ij}$'s are nonnegative numbers. The systems defined in this fashion are described by polynomial dynamics. The only assumptions required in order for the results to hold are that the $b_\ell$'s be linearly independent and that the $m \times m$ matrix $A = (a_{ij})$ must be *irreducible*. Recall that this means that $(I + A)^{m-1} \in \mathbb{R}_+^{m \times m}$ or, equivalently, that the incidence graph $G(A)$ is strongly connected (where $G(A)$ is the graph whose nodes are the integers $\{1, \ldots, m\}$ and for which there is an edge $j \to i$, $i \neq j$, if and only if $a_{ij} > 0$). This assumption amounts to a "weak reversibility" property in the application to chemical reactors, as discussed briefly in Section 3, and is crucial to the validity of the results.

We now describe the underlying dynamics of the systems to be studied, and leave for later the introduction of additional terms in order to model the possibility of control actions. Our systems are parametrized by two matrices $A$ and $B$ with nonnegative entries, as well as a collection of nonnegative functions $\theta_i$, $i = 1, \ldots, n$, and have the following general form:

$$\dot{x} \;=\; f(x) \;=\; \sum_{i=1}^{m}\sum_{j=1}^{m} a_{ij}\, \theta_1(x_1)^{b_{1j}} \theta_2(x_2)^{b_{2j}} \ldots \theta_n(x_n)^{b_{nj}}\, (b_i - b_j) \tag{5}$$

where $b_\ell$ denotes the $\ell$-th column of $B$ (notice that the diagonal entries of $A$ are irrelevant, since $b_i - b_i = 0$). Several restrictions on $A$, $B$, and the $\theta_i$'s are imposed below. The powers are interpreted as follows, for any $r, c \geq 0$: $r^0 = 1$, $0^c = 0$ if $c > 0$, and $r^c = e^{c \ln r}$ if $r > 0$ and $c > 0$.

The main motivating example, arising from mass-action kinetics in chemistry, is obtained when $\theta_i(y) = |y|$ for all $i$ and $B$ is a matrix whose entries are nonnegative integers (so, for nonnegative vectors $x$, we have polynomial equations). This is the case mentioned in the first paragraph, and will be referred to as the "standard setup" in this paper.

The hypotheses on the $\theta_i$'s, $A$, and $B$ are as follows. Each map

$$\theta_i \;:\; \mathbb{R} \to [0, \infty)$$

is locally Lipschitz, has $\theta_i(0) = 0$, satisfies $\int_0^1 |\ln \theta_i(y)|\, dy < \infty$, and its restriction to $\mathbb{R}_{\geq 0}$ is strictly increasing and onto. We suppose that

$$A \;=\; (a_{ij}) \in \mathbb{R}_{\geq 0}^{m \times m} \quad \text{is irreducible} \tag{6}$$

and, for $B$, that:

$$\text{each entry of } B \text{ is either } 0 \text{ or } \geq 1 \tag{7}$$

(this last hypothesis insures that $f(x)$ in (5) is a locally Lipschitz vector field, so we have uniqueness of solutions for the differential equation),

$$B \;=\; (b_1, \ldots, b_m) \in \mathbb{R}_{\geq 0}^{n \times m} \quad \text{has rank } m \tag{8}$$

(so, its columns $b_i$ are linearly independent), and

$$\text{no row of } B \text{ vanishes.} \tag{9}$$

This last hypothesis is made mainly for convenience. Observe that if some row, let us say the $k$-th one, were zero, then the same dynamics would be obtained if all entries in row $k$ are replaced by "1" (since the differences $b_{ki} - b_{kj}$ are still zero) and one restricts the dynamics to those states satisfying $x_k \equiv c$, where $c$ is the positive number such that $\theta_k(c) = 1$.

*From now on, we assume that all systems (5) considered satisfy the above assumptions.*

Our study will focus on those solutions of (5) which evolve in the nonnegative orthant $\mathbb{R}^n_{\geq 0}$. Recall that a subset $S \subseteq \mathbb{R}^n$ is said to be *forward invariant* with respect to the differential equation $\dot{x} = f(x)$ provided that each solution $x(\cdot)$ with $x(0) \in S$ has the property that $x(t) \in S$ for all positive $t$ in the domain of definition of $x(\cdot)$. It is routine to show (cf. Section 7) that the nonnegative and positive orthants are forward invariant:

**Lemma 2.1** Both $\mathbb{R}^n_{\geq 0}$ and $\mathbb{R}^n_+$ are forward-invariant sets with respect to the system (5).

(These properties are simple consequences of the fact that, because of the assumptions made, the $k$-th component of a solution of (5) will satisfy $\dot{x}_k(t) \geq 0$ whenever $x_k(t) = 0$.) We will also show in Section 8.2.4 that there are no finite explosion times:

**Lemma 2.2** For each $\xi \in \mathbb{R}^n_{\geq 0}$ there is a (unique) solution $x(\cdot)$ of (5) with $x(0) = \xi$, defined for all $t \geq 0$.

In order to state concisely the main results for systems (5), we need to introduce a few additional objects. The subspace

$$\mathcal{D} := \text{span}\,\{b_i - b_j,\, i \neq j\} \;=\; \text{span}\,\{b_1 - b_2, \ldots, b_1 - b_m\} \tag{10}$$

can be seen as a distribution in the tangent space of $\mathbb{R}^n$; it has dimension $m-1$ because adding $b_1$ to the last-shown generating set gives the column space of the rank-$m$ matrix $B$. For each vector $p \in \mathbb{R}^n$, we may also consider the parallel translate of $\mathcal{D}$ that passes through $p$, i.e. $p + \mathcal{D} = \{p + d,\, d \in \mathcal{D}\}$. A set $S$ which arises as an intersection of such an affine subspace with the nonnegative orthant:

$$S \;=\; (p + \mathcal{D}) \bigcap \mathbb{R}^n_{\geq 0}$$

(for some $p$, without loss of generality in $\mathbb{R}^n_{\geq 0}$) will be referred to as a *class*. If $S$ intersects the positive orthant $\mathbb{R}^n_+$, we say that $S$ is a *positive class*. The significance of classes is given by the fact that any solution $x(\cdot)$ of (5) must satisfy

$$x(t) - x(0) \;=\; \int_0^t f(x(s))\, ds \;=\; \sum_{i=1}^m \sum_{j=1}^m \kappa(t)\,(b_i - b_j) \;\in \mathcal{D}\,, \tag{11}$$

$\kappa(t) = \int_0^t a_{ij} \theta_1(x_1(s))^{b_{1j}} \theta_2(x_2(s))^{b_{2j}} \ldots \theta_n(x_n(s))^{b_{nj}}\, ds$, so $x(t) \in x(0) + \mathcal{D}$ for all $t$. In particular:

**Lemma 2.3** Each class is forward invariant for (5). □

The introduction of control action, through additional feedback loops, may be used in order to overcome the constraints imposed by (11). In order to formulate our control-theoretic results, we will suppose that $r$ external inputs $u_\ell$ can be used to independently influence each of $\ell$ state coordinates, In other words, we will also consider the following control system associated to the basic open-loop model (5):

$$\dot{x} \;=\; f(x) + \sum_{\ell=1}^r u_\ell\, e_{k_\ell} \tag{12}$$

($f$ is as in (5)), where $r$ is a positive integer, $e_1, \ldots, e_n$ are the $n$ canonical basis vectors in $\mathbb{R}^n$, and $k_1, \ldots, k_r$ are $r$ distinct elements of $\{1, \ldots, n\}$.

Of course, inputs might influence the system in manners other than through independent action on some coordinates, which represent what are sometimes (especially in compartmental models) called *inflow-controlled systems*, see [1]. Bilinear control action, for instance, is arguably more interesting in reaction systems. We will only study the above class of control systems, but research is ongoing on generalizations to other formulations.

We will first state the main results for autonomous systems (5), and later we state a stabilization result for (12).

We denote by $E$ (respectively, $E_+$ or $E_0$) the set of nonnegative (respectively, positive or boundary) equilibria of (5), i.e. the set of states $\bar{x} \in \mathbb{R}^n_{\geq 0}$ (respectively, $\in \mathbb{R}^n_+$ or $\in \mathbb{R}^n_0$), such that $f(\bar{x}) = 0$. Of course, $E$ is the disjoint union of $E_+$ and $E_0$.

**Theorem 1** *Consider any system (5), under the stated assumptions. For every maximal solution of (5) with $x(0) \in \mathbb{R}^n_{\geq 0}$, it holds that $x(t) \to E$ as $t \to +\infty$.*

This will be proved in Section 8.2.2. The invariance of classes (which are contained in subspaces of dimension $m - 1 < n$) precludes asymptotic stability of equilibria of (5). The appropriate concept is that of *asymptotic stability relative to a class*. We say that an equilibrium $\bar{x} \in S$ is asymptotically stable relative to a class $S$ if it is (a) stable relative to $S$ (for each $\varepsilon > 0$, there is some $\delta > 0$ such that, for all solutions $x(\cdot)$, $|x(0) - \bar{x}| < \delta$ and $x(0) \in S$ imply $|x(t) - \bar{x}| < \varepsilon$ for all $t \geq 0$) and (b) locally attractive relative to $S$ (for some $\varepsilon > 0$, if $|x(0) - \bar{x}| < \varepsilon$ and $x(0) \in S$ then $x(t) \to \bar{x}$ as $t \to +\infty$). We say that $\bar{x} \in S$ is *globally* asymptotically stable relative to a class $S$ if it is stable relative to $S$ and globally attractive relative to $S$ ($x(t) \to \bar{x}$ for all solutions with $x(0) \in S$). The main results are as follows; the first part is shown in Section 5.1, and the remaining two in Section 8.

**Theorem 2** *Consider any system (5), under the stated assumptions. Fix any positive class $S$.*

a. *There is a unique equilibrium $\bar{x}_S \in S \cap E_+$.*

b. *The equilibrium $\bar{x}_S$ is asymptotically stable relative to $S$.*

c. *The equilibrium $\bar{x}_S$ is globally asymptotically stable relative to $S$ if and only if $S \cap E_0 = \emptyset$.*

**Example 2.4** The following two trivial examples may help in understanding the above theorems. In both cases we take $n = m = 2$ and $\theta_i(y) = |y|$ for $i = 1, 2$. The first example has

$$A = \begin{pmatrix} 0 & 1 \\ 1 & 0 \end{pmatrix}, \quad B = \begin{pmatrix} 1 & 2 \\ 1 & 1 \end{pmatrix}.$$

The system (5) is (for nonnegative states):

$$\dot{x}_1 = (x_1 - 1)x_1 x_2$$
$$\dot{x}_2 = 0$$

and thus $E_0 = \mathbb{R}^2_0 = \{x \mid x_1 x_2 = 0\}$ and $E_+ = \{x \mid x_1 = 1, x_2 > 0\}$. The positive classes are the sets $S = S_r = \{x \mid x_1 \geq 0, x_2 = r\}$, for each $r > 0$, and for each such $S = S_r$, $\bar{x}_S = (1, r)'$ is asymptotically stable with domain of attraction $\{x \mid x_1 > 0, x_2 = r\}$. See Figure 1. Each class $S_r$ has a second equilibrium $(0, r)'$, but this second equilibrium is in the boundary, so there is no contradiction with part $a$ of Theorem 2. Regarding Theorem 1, observe that every trajectory either converges to an interior equilibrium $(1, r)'$ or it is itself a trajectory consisting of an equilibrium (and hence also converges to $E$, in a trivial sense).

The second example has

$$A = \begin{pmatrix} 0 & 1 \\ 1 & 0 \end{pmatrix}, \quad B = \begin{pmatrix} 1 & 0 \\ 0 & 1 \end{pmatrix}.$$

We now have (for nonnegative states) the following linear system:

$$\dot{x}_1 = x_2 - x_1$$
$$\dot{x}_2 = x_1 - x_2$$

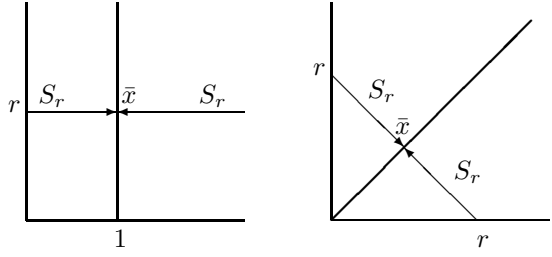

Figure 1: Equilibria for example 2.4, dark lines indicate equilibria

and thus $E_0 = \{0\}$ and $E_+ = \{x \,|\, x_1 = x_2 > 0\}$. The positive classes contain no boundary equilibria; they are the sets $S = S_r = \{x \in \mathbb{R}^2_{\geq 0} \,|\, x_1 - x_2 = r\}$, for each $r > 0$ (see Figure 1). For each such $S = S_r$, $\bar{x}_S = (r, r)'$ is globally asymptotically stable. □

**Feedback stabilization**

The invariance of classes precludes the existence of *global* (not merely relative to a class) attractors for the uncontrolled system (5). Now suppose that we wish to find a feedback law which forces all the solutions of the closed-loop system obtained from the controlled system (12) under this feedback to converge, as $t \to +\infty$, to a specific equilibrium $\bar{x} \in \mathbb{R}^n_+$. It turns out that the obvious negative feedback solution, namely to use inputs proportional to the errors $\bar{x}_k - x_k$, achieves the goal of global stabilization to $\bar{x}$, provided that enough inputs are used ($r \geq n - m + 1$, which is obviously necessary, since solutions of (12) evolve in linear subspaces of dimension $m - 1 + r$) and that the $k_\ell$'s are appropriately chosen. We next state a result in that regard. For each $j \in \{1, \ldots, m\}$ we consider the set

$$\mathcal{S}_j := \{k \,|\, b_{kj} > 0\} \tag{13}$$

which is nonempty, by (8).

**Theorem 3** *Let $r = n - m + 1$, and suppose that $\gamma_1, \ldots, \gamma_r$ are arbitrary positive real numbers, and $k_1, \ldots, k_r \in \{1, \ldots, n\}$ are such that*

$$\mathcal{D} + \mathrm{span}\,\{e_{k_1}, \ldots, e_{k_r}\} = \mathbb{R}^n \tag{14}$$

*and*

$$\mathcal{S}_j \subseteq \{k_1, \ldots, k_r\} \tag{15}$$

*for some $j \in \{1, \ldots, m\}$. Pick any equilibrium $\bar{x} \in E_+$. Then, all maximal solutions of*

$$\dot{x} = f(x) + \sum_{\ell=1}^{r} \gamma_\ell (\bar{x}_{k_\ell} - x_{k_\ell}) \, e_{k_\ell} \tag{16}$$

*with $x(0) \in \mathbb{R}^n_{\geq 0}$ are defined for all $t \geq 0$ and remain in $\mathbb{R}^n_{\geq 0}$, and $\bar{x}$ is a globally asymptotically stable equilibrium of (16).*

The invariance statement is proved in Section 7. The global stability statement is proved in Section 8.2.1.

To study the closed-loop system (16), and also to be able to formulate a result concerning robustness of the stability properties described by Theorems 1 and 2, we will study "positive perturbations" of the basic uncontrolled system model (5). These are described by equations as follows:

$$\dot{x} = f^*(x) = f(x) + g(x) \tag{17}$$

where $f$ is as in (5) (and all the stated assumptions hold), and $g$ is a locally Lipschitz vector field on $\mathbb{R}^n$ for which the following property holds:

$$(\forall x \in \mathbb{R}^n_{\geq 0}) \ (\forall k \in \{1, \ldots, n\}) \ [x_k = 0 \Rightarrow g_k(x) \geq 0] \tag{18}$$

where $g_k(x)$ denotes the $k$-th coordinate of $g(x)$. Property (18) is the most natural assumption which guarantees the forward-invariance of $\mathbb{R}^n_{\geq 0}$ (since it will imply that $\dot{x}_k(t) \geq 0$ whenever $x_k(t) = 0$). Of course, any results established for arbitrary systems (17) will be also true for systems (5) (take $g \equiv 0$).

The feedback system (16) is of this type. To see that (18) holds, note that either $g_k(x) = \gamma_\ell(\bar{x}_{k_\ell} - x_{k_\ell})$ (if $k = $ some $k_\ell$) or $g_k(x) = 0$, and, in the first case, $x_k = 0$ implies $g_k(x) = \gamma_\ell \bar{x}_{k_\ell} > 0$. Moreover, it satisfies the following strengthening of (18):

$$(\exists j \in \{1, \ldots, m\}) \ (\forall x \in \mathbb{R}^n_{\geq 0}) \ (\forall k \in \mathcal{S}_j) \ [x_k = 0 \Rightarrow g_k(x) > 0] \tag{19}$$

(note the strict inequality, in contrast to (18)). Indeed, pick $j$ as in (15), $k \in \mathcal{S}_j$, and $x$ so that $x_k = 0$. Then, for each $k \in \mathcal{S}_j$, property (15) guarantees that $g_k(x) = \gamma_\ell \bar{x}_{k_\ell} > 0$.

**Robustness**

A different specialization of the general form (17) allows the study of robustness with respect to perturbations which preserve classes. The corresponding systems are obtained by adding vector fields which lie pointwise in the span of the $b_i - b_j$'s. We suppose given a collection of locally Lipschitz functions $\Delta_{ij} : \mathbb{R}^n \to \mathbb{R}_{\geq 0}$ $(i, j \in \{1, \ldots, m\})$ such that

$$(\forall i, j \in \{1, \ldots, m\}) \ (\forall x \in \mathbb{R}^n_{\geq 0}) \ (\forall k \in \mathcal{S}_j) \ [x_k = 0 \Rightarrow \Delta_{ij}(x) = 0] \tag{20}$$

and, using these, define the system:

$$\dot{x} = f(x) + \sum_{i=1}^{m} \sum_{j=1}^{m} \Delta_{ij}(x) (b_i - b_j) . \tag{21}$$

Observe that property (20) implies that this system is of the general type (17). Indeed, the only possible negative signs for $g_k(x) = \sum \sum \Delta_{ij}(x)(b_{ki} - b_{kj})$ can arise from the terms of the form $\Delta_{ij}(x)b_{kj}$ with $b_{kj} \neq 0$ (i.e. $k \in \mathcal{S}_j$), but these vanish, because of (20), when $x_k = 0$. Of course, systems (5) are a subclass of (21) (take all $\Delta_{ij} \equiv 0$). The main robustness result will be as follows.

**Theorem 4** *For each positive class $S$ there exists a continuous function $\delta_S : S \bigcap \mathbb{R}^n_+ \to \mathbb{R}_{\geq 0}$, with $\delta_S(x) > 0$ if and only if $x \neq \bar{x}_S$, such that, for any collection $\{\Delta_{ij}\}$ such that*

$$\sum_{i=1}^{m} \sum_{j=1}^{m} \Delta_{ij}(x)^2 \leq \delta_S(x) \tag{22}$$

*for all $x \in S \bigcap \mathbb{R}^n_+$, the following properties hold for the system (21):*

1. *Both $\mathbb{R}^n_{\geq 0}$ and $\mathbb{R}^n_+$, and the class $S$, are forward-invariant.*

2. *For each $\xi \in S$ there is a (unique) solution $x(\cdot)$ with $x(0) = \xi$, defined for all $t \geq 0$.*

3. *The equilibrium $\bar{x}_S$ is asymptotically stable relative to $S$.*

4. *The equilibrium $\bar{x}_S$ is globally asymptotically stable relative to $S$ if and only if $S \bigcap E_0 = \emptyset$.*

In other words, the claim is that the conclusions previously stated for systems (5) are preserved by perturbations, as long as the magnitude constraint (22) is satisfied. The result is not in itself surprising; the main interest is in its proof, which will be largely constructive, offering an explicit formula for the function $\delta_S$ and, moreover, the fact that $\delta_S$ will depend nicely on the class $S$ (in a sense that will be clear). The invariance of $S$ is clear from the form (21), and the (also easy) invariance of orthants is proved in Section 7. Parts 2-4 are proved in Section 8.2.3.

**Regularity**

For all results to be stated next, we suppose a fixed system (5) has been given.

We say that a function defined on an open subset of $\mathbb{R}^n$ is of class $\mathcal{C}^k$, $k = 1, 2, \ldots, \infty, \omega$ if it is $k$-times continuously differentiable (for $k = 1, \ldots, \infty$) or real-analytic (for $k = \omega$). We will refer to the following hypothesis, for each such $k$:

$$(\text{H}_k) \quad \text{each } \theta_i \text{ restricted to } \mathbb{R}_+ \text{ is of class } \mathcal{C}^k, \text{ and } \theta'_i(y) > 0 \text{ for all } y > 0.$$

Obviously, this hypothesis is satisfied with $k = \omega$ when $\theta_i$ is the identity map for $y > 0$, as is the case in the standard setup leading to polynomial systems.

We will show, in Section 5.2:

**Theorem 5** *If hypothesis ($\text{H}_k$) holds, then $E_+$ is an embedded submanifold of $\mathbb{R}_+^n$, $\mathcal{C}^k$-diffeomorphic to $\mathbb{R}^{n-m+1}$.*

Observe that, in the standard setup, $f$ is a polynomial vector field, so $E_+ = \{x \in \mathbb{R}_+^n \mid f(x) = 0\}$ and hence is an algebraic subset.

Theorem 2 states that there is a unique interior equilibrium $\bar{x}$ in each positive class. This equilibrium depends nicely (smoothly, analytically) on the class, provided that the $\theta_i$'s be regular enough. In order to make this statement precise, we introduce the map $\pi : \mathbb{R}_+^n \to \mathbb{R}_+^n$ which assigns, to each $x \in \mathbb{R}_+^n$, the unique interior equilibrium $\bar{x}_S$ in the class $S$ which contains $x$. The next result is proved in Section 5.1.

**Theorem 6** *If hypothesis ($\text{H}_k$) holds, then $\pi$ is of class $\mathcal{C}^k$.*

Note that, once that we know that $E_+$ an embedded submanifold, there is no ambiguity in the above statement: if we view $\pi$ as a map $\pi : \mathbb{R}_+^n \to E_+$, it is also of class $\mathcal{C}^k$.

The following open subset of $\mathbb{R}^n$:

$$\mathcal{R} := \mathbb{R}_+^n + \mathcal{D} \qquad (23)$$

is the union of all those parallel translates of $\mathcal{D}$ which intersect the positive orthant. It includes all positive classes. Assuming again hypothesis ($\text{H}_k$), it turns out that one may always find a change of variables which transforms $\mathcal{R}$ onto $\mathbb{R}^n$, and in particular $E_+$ onto the set $\{(X_1, X_2) \mid X_1 = 0\}$ and positive classes into subsets of sets of the form $\{(X_1, X_2) \mid X_2 = c\}$ for constants $c$, and transforms the dynamics on $\mathcal{R}$ into the following form:

$$\dot{X}_1 = F_1(X_1, X_2)$$
$$\dot{X}_2 = 0$$

where $X_1$ and $X_2$ represent blocks of variables according to the decomposition $\mathbb{R}^n = \mathbb{R}^{m-1} \times \mathbb{R}^{n-m+1}$. Since points of $E_+$ are equilibria, this implies that $F_1(0, X_2) \equiv 0$. To be precise, we prove in Section 5.3:

**Theorem 7** *Assume that hypothesis ($\text{H}_k$) holds. Then, there exists a $\mathcal{C}^k$ diffeomorphism $\Phi : \mathcal{R} \to \mathbb{R}^{m-1} \times \mathbb{R}^{n-m+1}$ such that, denoting*

$$F(X) = \Phi_*(\Phi^{-1}(X)) f(\Phi^{-1}(X))$$

*and writing $F = (F_1, F_2)$ and $\Phi = (\Phi_1, \Phi_2)$ in block form,*

1. *$x - z \in \mathcal{D}$ if and only if $\Phi_2(x) = \Phi_2(z)$,*

2. *$\Phi_1(x) = 0$ if and only if $x \in E_+$, and hence $F(0, X_2) = 0$ for all $X_2$, and*

3. $F_2(X) = 0$ *for all $X$.*

The preceding results provide steps in proving this one, but conversely, Theorem 7 implies Theorem 5, since $E_+$ is transformed into a coordinate space. In fact, it also establishes that $E_+$ transversally intersects each class $S$. Also Theorem 7 implies Theorem 6, since the map that selects the element in $E_+$ in the same class as a given point amounts to a coordinate projection, let us call it $\widetilde{\pi}$, onto the last $n - m + 1$ variables under the change of variables given by $\Phi$. Moreover, let $\mathcal{O}$ be the set of nonnegative points that are not boundary equilibria of (5), that is, $\mathcal{O} = \mathbb{R}^n_{\geq 0} \setminus E_0$. We will prove (Corollary 7.7) that each $x \in \mathcal{O}$ belongs to some positive class; thus, the mapping $\pi$ extends to $\mathcal{O}$. Since $\mathcal{R}$ includes every positive class, it includes $\mathcal{O}$, and the extension of $\pi$ to $\mathcal{O}$ transforms into the restriction of $\widetilde{\pi}$ to $\Phi(\mathcal{O})$. Thus, we have the following consequence of Theorem 7:

**Corollary 2.5** *If hypothesis ($H_k$) holds, then the extension of $\pi$ to $\mathcal{O}$ is also of class $\mathcal{C}^k$.* □

## 3 Mass-Action Kinetics, and McKeithan's System

As mentioned in the Introduction, the results explained in Theorems 1 and 2 are basically theorems for what are called mass-action networks of zero deficiency, and are given (implicitly in the case of global stability) in [9, 10, 11, 12, 16]. (There is one global stability result stated explicitly in the above papers, namely in [16]. The statement would be that every trajectory which starts in the positive orthant must converge to the interior equilibrium in the corresponding class. However, this question is still open, since the suggested proof used the implication that if a positive definite function $V$ of states has a negative derivative along trajectories while away from an equilibrium, then global stability of the equilibrium follows. This implication is, in general, false, without the assumption of some sort of radial unboundedness, a property which definitely does not hold in the context in which it is being applied; see [17] for a retraction of that proof.)

The assumptions of irreducibility of $A$ and full rank of $B$ (both of which may be relaxed somewhat, cf. Remark 5.5) are key ones. They serve to rule out periodic (or even chaotic) behaviors which may otherwise arise in chemical networks such as the Belousov-Zhabotinsky reaction or Prigogine-Lefever's "Brusselator" (for which see e.g. [7]).

### 3.1 Kinetic Proofreading in T-Cell Signaling

The equations with which we started represent the dynamics of the "kinetic proofreading" model proposed by McKeithan in [19] in order to describe how a chain of modifications of the T-cell receptor complex, via tyrosine phosphorylation and other reactions, may give rise to both increased sensitivity and selectivity of response. Let us introduce two additional variables $T(t)$ and $M(t)$, which represent the concentrations of T-cell receptor (TCR) and a peptide-major histocompatibility complex (MHC). The constant $k_1$ is the association rate constant for the reaction which produces an initial ligand-receptor complex $C_0$ from TCR's and MHC's. The quantities $C_i(t)$ represent concentrations of various intermediate complexes, and McKeithan postulates that recognition signals are determined by the concentrations of the final complex $C_N$. The constants $k_{p,i}$ are the rate constants for each of the steps of phosphorylation or other intermediate modifications, and the constants $k_{-1,i}$ are dissociation rates.

Global stability of a unique equilibrium will be deduced from Part c in Theorem 2 when we view the equations as those of an appropriate system (5) restricted to a suitable class (which is determined by the constants $M^*$ and $T^*$). The complete reaction network is represented graphically in Figure 2 (this is the same as Figure 1 in [19], except that we do not make the simplifying assumption of equal rates). In other words, we have a system of form (5), where we use $x = (T, M, C_0, \ldots, C_N)'$ as a state, and take

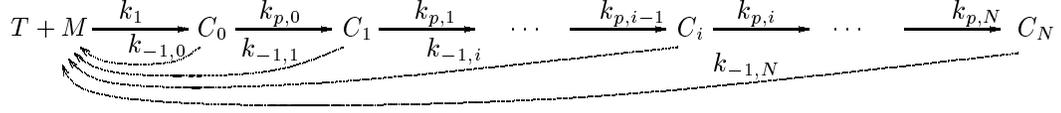

Figure 2: McKeithan's network

$m = n - 1 = N + 2$,

$$b_1 = \begin{pmatrix} 1 \\ 1 \\ 0 \\ 0 \\ \vdots \\ 0 \end{pmatrix}, \ b_2 = \begin{pmatrix} 0 \\ 0 \\ 1 \\ 0 \\ \vdots \\ 0 \end{pmatrix}, \ b_3 = \begin{pmatrix} 0 \\ 0 \\ 0 \\ 1 \\ \vdots \\ 0 \end{pmatrix}, \ \ldots \ b_m = \begin{pmatrix} 0 \\ 0 \\ 0 \\ 0 \\ \vdots \\ 1 \end{pmatrix},$$

and $A = (a_{ij})$ with $a_{21} = k_1$, $a_{1i} = k_{-1,i-2}$ ($i = 2, \ldots, m$), $a_{i,i-1} = k_{p,i-3}$ ($i = 3, \ldots, m$), and all other $a_{ij} = 0$. The corresponding set of differential equations is:

$$\begin{aligned}
\dot{T} &= -k_1 TM + \sum_{i=0}^{N} k_{-1,i} C_i \\
\dot{M} &= -k_1 TM + \sum_{i=0}^{N} k_{-1,i} C_i \\
\dot{C}_0 &= k_1 TM - (k_{-1,0} + k_{p,0}) C_0 \\
&\vdots \\
\dot{C}_i &= k_{p,i-1} C_{i-1} - (k_{-1,i} + k_{p,i}) C_i \\
&\vdots \\
\dot{C}_N &= k_{p,N-1} C_{N-1} - k_{-1,N} C_N .
\end{aligned}$$

Notice that $\mathcal{D} = \{x \mid T + C_0 + \ldots + C_N = M + C_0 + \ldots + C_N = 0\}$, and the positive classes are of the form $S = S^{\alpha,\beta}$, intersections with $\mathbb{R}^n_{\geq 0}$ of the affine planes

$$T + C_0 + \ldots + C_N = \alpha, \quad M + C_0 + \ldots + C_N = \beta,$$

with $\alpha > 0$ and $\beta > 0$. The original system is nothing else than the restriction of the dynamics to the class determined by $\alpha = T^*$ and $\beta = M^*$. Thus, the conclusions will follow from Theorem 2 as soon as we prove that $S \bigcap E_0 = \emptyset$ for any positive class $S$. Pick any $x \in \mathbb{R}^n_0$ and any positive class $S^{\alpha,\beta}$. According to Proposition 6.3, it will be enough to find some $j \in \{1, \ldots, m\}$ with the property that $x_k \neq 0$ for all $k \in \mathcal{S}_j$. Here, $\mathcal{S}_1 = \{1, 2\}$ and $\mathcal{S}_j = \{j + 1\}$ for $j = 2, \ldots, m$. If the property is not satisfied for some $j \in \{2, \ldots, m\}$, then $C_i = 0$ for all $i$. But in this case, the equations for $S^{\alpha,\beta}$ give that $T = \alpha > 0$ and also $M = \beta > 0$, so $j = 1$ can be used. In conclusion, $x \notin E_0$, and hence Part c of the theorem applies.

Let us compute equilibria explicitly. Setting right-hand sides to zero gives $C_N = (k_{p,N-1}/k_{-1,N})C_{N-1}$, and recursively using $C_i = (k_{p,i-1}/(k_{-1,i} + k_{p,i}))C_{i-1}$ we may express all $C_i$'s as multiples of $C_0$. We also have $C_0 = [k_1/(k_{-1,0} + k_{p,0})]TM$. Thus

$$E_+ = \{(\alpha, \beta, \kappa_0 \alpha \beta, \kappa_1 \alpha \beta, \ldots, \kappa_N \alpha \beta) \mid \alpha, \beta > 0\}$$

where the $\kappa_i$'s are rational functions of the constants defining the system (these are the equilibria studied in [19]). It is obvious in this example that we obtained a two-dimensional ($n - m + 1 = 2$) nonsingular algebraic subvariety of $\mathbb{R}^n$.

We next discuss the application of Theorem 3 to this example. As $r = 2$, we consider feedback laws of the type $\gamma_1(\bar{x}_{k_1} - x_{k_1}) + \gamma_2(\bar{x}_{k_2} - x_{k_2})$, where $\gamma_i > 0$ and $k_1$, $k_2$ satisfy (14) and (15). Since $\mathcal{D}$ is spanned by $e_1 + e_2 - e_3, \ldots, e_1 + e_2 - e_n$, this means that (14) is satisfied provided that $1 \in \{k_1, k_2\}$ or $2 \in \{k_1, k_2\}$, and this is sufficient to guarantee that also (15) satisfied. In other words, stabilization to any desired equilibrium is possible provided that some pair of reactants, including at least one of $T$ or $M$, can be manipulated by a controller.

To illustrate some of our constructions in a comparatively trivial case, let us now specialize even further, taking just $N = 0$ and all constants equal to one. (This is mathematically the same as the example treated in Section 1.1.) Writing $x, y, z$ instead of $T, M, C_0$, we have the system

$$\dot{x} = -xy + z \,,\; \dot{y} = -xy + z \,,\; \dot{z} = xy - z$$

and $E_+ = \{(x, y, z) \,|\, z - xy = 0\} \bigcap \mathbb{R}^3_+$ is a hyperbolic paraboloid (intersected with the main orthant). This is obviously a nonsingular algebraic set, but it is instructive to see it also as the image of the diffeomorphism $M : \mathcal{D}^\perp \to E_+$ constructed in the proof of Theorem 5 in Section 5.2. Using $\bar{x} = (1, 1, 1)'$, we have the formula $M(r_1, r_2, r_3) = (e^{r_1}, e^{r_2}, e^{r_3})'$ for $(r_1, r_2, r_3) \in \mathcal{D}^\perp$. Here, $\mathcal{D}$ is the span of $(1, 1, -1)'$, i.e. $\{(r_1, r_2, r_3) \,|\, r_1 + r_3 = r_2 + r_3 = 0\}$, and $\mathcal{D}^\perp = \{(r_1, r_2, r_3) \,|\, r_3 = r_1 + r_2\}$, so it is clear that $M(\mathcal{D}^\perp) = E_+$.

Note, incidentally. that there are boundary equilibria as well: $E_0 = \{(x, 0, 0) \,|\, x \geq 0\} \bigcup \{(0, y, 0) \,|\, y \geq 0\}$, but none of these is in a positive class (since both $x + z$ and $y + z$ are positive constants on positive classes).

The most natural (in the motivating application, anyway) initial state is one in which $z(0) = 0$. We now compute the value of the map $\pi$ considered in Corollary 2.5 at such a point $(0, y_0, z_0)$. We must find the positive-orthant intersection of $\{z = xy\}$ with the line $L_{x_0, y_0} = \{x + z = x_0,\, y + z = y_0\}$. Equivalently, we need to solve $(z - x_0)(z - y_0) - z = 0$ subject to the constraint $0 < z < \min\{x_0, y_0\}$, which guarantees that $x > 0$ and $y > 0$ as well as $z > 0$. There is exactly one such solution, and it is the smallest of the two solutions, since the graphs of $f_1(z) = z$ and $f_2(z) = (z - x_0)(z - y_0)$ intersect at precisely one point in the interval $(0, \min\{x_0, y_0\})$, so we take the negative sign in the quadratic formula:

$$\pi_3(x_0, y_0, 0) = \frac{1}{2}\left[x_0 + y_0 + 1 - \sqrt{(x_0 - y_0)^2 + 2(x_0 + y_0) + 1}\right]$$

(and corresponding values for $\pi_1$ and $\pi_2$), which is, indeed, a real-analytic (hypothesis H$_\omega$) function.

Finally, let us find a decomposition as in Theorem 7. We let $X_1 := z$ and $X_2 = (u, v) = (x + z, y + z)$. Thus $\dot{X}_2 \equiv 0$ along all solutions, and $\dot{X}_1 = \dot{z} = F_1(z, u, v) = (u - z)(v - z) - z$. This does not yet have $F_1(0, u, v) = 0$ for all $u, v$, so one needs a further change of variables to bring the equilibrium to zero, which can be achieved by translating $z \mapsto \tilde{z} := z - \pi_3(0, u, v)$. Note that $E_+$, the set of points with $z = xy$, gets mapped into $\tilde{z} = 0$, and the positive orthant gets mapped into the set given by the inequalities $u > 0$, $v > 0$, and $z \in (-\pi_3(0, u, v), \min\{u, v\} - \pi_3(0, u, v))$.

## 4 Some Preliminary Facts

The equations (5) (and their perturbed form (17)) have a considerable amount of structure, and various useful properties are reflected in alternative expressions for the system equations.

### 4.1 Other Expressions for the System Equations

For each $i = 1, \ldots, n$, let us introduce the maps

$$\rho_i(y) \;:=\; \ln \theta_i(y)$$

(with $\rho_i(0) = -\infty$); note that $\lim_{y \searrow 0} \rho_i(y) = -\infty$ and $\int_0^1 |\rho_i(y)|\, dy < \infty$. Furthermore, the restriction of $\rho_i$ to $\mathbb{R}_+$ is locally Lipschitz, strictly increasing, and onto $\mathbb{R}$. For any positive integer $n$, we let

$$\vec{\rho}\,:\, \mathbb{R}^n \to [-\infty, \infty)^n \,:\, x \mapsto (\rho_1(x_1), \ldots, \rho_n(x_n))'$$

(we do not write "$\vec{\rho}_n$" to emphasize the dependence on $n$, because $n$ will be clear from the context). Then (5) can also be written as

$$\dot{x} \;=\; f(x) \;=\; \sum_{i=1}^m \sum_{j=1}^m a_{ij} e^{\langle b_j, \vec{\rho}(x)\rangle} (b_i - b_j)\,. \tag{24}$$

The expression "$e^{\langle b_j, \vec{\rho}(x)\rangle}$" is interpreted in accordance with the conventions made for powers: if $x$ is a vector and $k \in \{1, \ldots, n\}$ is an index such that $x_k = 0$ and $b_{kj} > 0$, then $e^{b_{kj} \rho_k(x_k)} = 0$, consistently with $e^{-\infty} = 0$, and thus also

$$e^{\langle b_j, \vec{\rho}(x)\rangle} \;=\; e^{b_{1j}\rho_1(x_1)} e^{b_{2j}\rho_2(x_2)} \ldots e^{b_{nj}\rho_n(x_n)} \;=\; 0\,,$$

but, if $b_{kj} = 0$, then we have $e^{b_{kj}\rho_k(x_k)} = 1$.

Another useful way of rewriting (5) is as follows. We write $f_k$ for the $k$-th coordinate of $f$ (i.e., the coordinates $x_k$ of solutions $x$ satisfy $\dot{x}_k = f_k(x)$). The terms in the sums defining $f_k$ can be collected into two disjoint sets: those that do not involve a product containing $\theta_k(x_k)$, for which $b_{kj} = 0$, and those which do involve $\theta_k(x_k)$. The latter, by assumption (7), have $b_{kj} \geq 1$, so we can factor $\theta_k(x_k)$ from $\theta_1(x_1)^{b_{1j}} \theta_2(x_2)^{b_{2j}} \ldots \theta_n(x_n)^{b_{nj}}$ and there remains a locally Lipschitz product. In other words, we can introduce, for each $k \in \{1, \ldots, n\}$, these two locally Lipschitz functions:

$$\alpha_k(x) \;:=\; \sum_{j \in J_{k,1}} \left( \sum_{i=1}^m a_{ij}\,(b_{ki} - b_{kj}) \right) \theta_1(x_1)^{b_{1j}} \theta_2(x_2)^{b_{2j}} \ldots \theta_k(x_k)^{b_{kj}-1} \ldots \theta_n(x_n)^{b_{nj}} \tag{25}$$

and

$$\beta_k(x) \;:=\; \sum_{j \in J_{k,0}} \left( \sum_{i=1}^m a_{ij} b_{ki} \right) \theta_1(x_1)^{b_{1j}} \theta_2(x_2)^{b_{2j}} \ldots \theta_n(x_n)^{b_{nj}}\,, \tag{26}$$

where $J_{k,1} := \{j \,|\, b_{kj} \geq 1\}$ and $J_{k,0} := \{j \,|\, b_{kj} = 0\}$. (Note that, in terms of the sets introduced in (13), $k \in \mathcal{S}_j$ if and only if $j \in J_{k,1}$.) Note that $\beta_k(x) \geq 0$ for all $x$. More generally, for perturbed systems (17), we let

$$\widetilde{\beta}_k \;:=\; \beta_k + g_k\,.$$

In terms of these functions,

$$\begin{aligned} f_k(x) &\;=\; \alpha_k(x)\,\theta_k(x_k) + \beta_k(x)\,, \\ f_k^*(x) &\;=\; \alpha_k(x)\,\theta_k(x_k) + \widetilde{\beta}_k(x) \end{aligned}$$

(where $f_k^* = f_k + g_k$). If $x \in \mathbb{R}^n_{\geq 0}$ and $k$ are such that $x_k = 0$, then property (18) says that $g_k(x) \geq 0$. In particular, since $\theta_k(0) = 0$,

$$x_k = 0 \;\Rightarrow\; f_k(x) = \beta_k(x) \geq 0 \;\text{ and }\; f_k^*(x) = \widetilde{\beta}_k(x) \geq 0 \tag{27}$$

so the vector fields $f$ and $f^* = f + g$ always point towards the nonnegative orthant, on the boundary $\mathbb{R}^n_0$.

## 4.2  A Coordinatization Property

When we apply the following Lemma, we will always take $D = \mathcal{D}$, but we can state the result in more generality. Actually, we may get an even more general result by dropping the assumption that the maps $\theta_i$ are onto. We will only assume, for the next result, that each $i = 1, \ldots, n$, $\theta_i : \mathbb{R} \to [0, \sigma_i)$ (where $0 < \sigma_i \leq \infty$) is locally Lipschitz, has $\theta_i(0) = 0$, satisfies $\int_0^1 |\ln \theta_i(y)|\, dy < \infty$, and its restriction to $\mathbb{R}_{\geq 0}$ is strictly increasing and onto $[0, \sigma_i)$, but we are not now requiring $\sigma_i = +\infty$. The reason for this relaxation is to allow consideration of functions like $\theta(x) = \frac{x}{k+x}$ which arise in Michaelis-Menten kinetics. We will impose, instead, another condition, which we describe next. Let, as earlier, $\rho_i(y) := \ln \theta_i(y)$, and let $\bar\rho_i := \ln \sigma_i$ (infinity if $\sigma_i = \infty$). Notice that $\rho_i^{-1}$, seen as the inverse of the restriction of $\rho_i$ to $\mathbb{R}_+$, is a strictly increasing map from $(-\infty, \bar\rho_i)$ onto $\mathbb{R}_+$. Thus, for any given constant $p$, $\rho_i^{-1}(s) > p + 1$ for all $s \geq t_0$, some $t_0 \in (-\infty, \bar\rho_i)$, which implies that, for $L(t) = \int_a^t \rho_i^{-1}(s)\, ds - pt$, its derivative $(dL/dt)(t) > 1$ for all $t \in (t_0, \bar\rho_i)$. Under the general assumption $\bar\rho_i = \infty$, we have that, for any $p$:

$$\lim_{t \nearrow \bar\rho_i} \int_a^t \rho_i^{-1}(s)\, ds - pt \;=\; +\infty \tag{28}$$

for any finite $a < \bar\rho_i$. Instead of assuming $\bar\rho_i = \infty$, we will merely ask that (28) should hold. (Example: $\theta(x) = \frac{x}{k+x}$ gives $\rho^{-1}(s) = \frac{ke^s}{1-e^s}$, so $L(t) = -k\ln(1 - e^t) - pt \to +\infty$ as $t \to 0 = \bar\rho$.)

**Lemma 4.1** Let $D$ be any subspace of $\mathbb{R}^n$. For each $p, q$ in $\mathbb{R}_+^n$, there exists a unique $x = \varphi(p, q) \in \mathbb{R}_+^n$ such that:

$$x - p \in D \tag{29}$$

and

$$\vec\rho(x) - \vec\rho(q) \in D^\perp. \tag{30}$$

Furthermore, if hypothesis $(H_k)$ holds, then the map $p, q \mapsto \varphi(p, q)$ is also of class $\mathcal{C}^k$.

*Proof.* We fix $p, q$ as in the statement, and start by introducing the following mapping, for each $i \in \{1, \ldots, n\}$:

$$L_i(t) := \int_0^{t + \rho_i(q_i)} \rho_i^{-1}(s)\, ds \;-\; p_i t$$

defined for $t \in \mathcal{L}_i := (-\infty, \bar\rho_i - \rho_i(q_i))$. By assumption (28), $L_i(t)$ increases to infinity as $t \nearrow \bar\rho_i - \rho_i(q_i)$. Also, $L_i(t) \to +\infty$ as $t \to -\infty$, because $\rho_i^{-1}$ is nonnegative and $p_i > 0$. Thus, $L_i$ is proper, that is, $\{t \mid L_i(t) \leq v\}$ is compact for each $v$.

Now we take the (continuously differentiable) function

$$Q(y) := \sum_{i=1}^n L_i(y_i)$$

thought of as a function of $y \in \mathcal{L} := \prod_{i=1}^n \mathcal{L}_i$. Observe that $0$ is in the open set $\mathcal{L}$, since $\bar\rho_i - \rho_i(q_i) > 0$ of all $i$ by definition of $\bar\rho_i$. This function is also proper, because

$$\{y \in \mathcal{L} \mid Q(y) \leq w\} \;\subseteq\; \prod_{i=1}^n \{t \in \mathcal{L}_i \mid L_i(t) \leq w - (n-1)\ell\}$$

where $\ell$ is any common lower bound for the functions $L_i$. Restricted to $D^\perp \cap \mathcal{L}$, $Q$ is still proper, so it attains a minimum at some point $y \in D^\perp \cap \mathcal{L}$, which is a nonempty (since $0$ belongs to it) relatively open subset of $D^\perp$. In particular, $y$ must be a critical point of $Q$ restricted to $D^\perp \cap \mathcal{L}$, so $(\nabla Q(y))' \in (D^\perp)^\perp = D$, which means that

$$\rho_1^{-1}(y_1 + \rho_1(q_1)) - p_1, \ldots, \rho_n^{-1}(y_n + \rho_n(q_n)) - p_n) \in D. \tag{31}$$

Finally, pick $x \in \mathbb{R}_+^n$ such that $\vec\rho(x) = y + \vec\rho(q)$. Such an $x$ exists because we can solve the equations $\rho_i(x_i) = y_i + \rho_i(q_i)$ for each $i$, because $y_i < \bar\rho_i - \rho_i(q_i)$ since $y \in \mathcal{L}$. Then $\vec\rho(x) - \vec\rho(q) \in D^\perp$ by definition, and (31) gives also $x - p \in D$.

Finally, we show uniqueness. Suppose that there is a second $z \in \mathbb{R}^n_+$ so that $z - p \in D$ and $\vec{\rho}(z) - \vec{\rho}(q) \in D^\perp$. This implies that $x - z \in D$ and $\vec{\rho}(x) - \vec{\rho}(z) \in D^\perp$. Since each $\rho_i$ is an increasing function, we have that, for any two distinct numbers $a, b$, $(a-b)(\rho_i(a) - \rho_i(b)) > 0$. So

$$\sum_{i=1}^n (x_i - z_i)(\rho_i(x_i) - \rho_i(z_i)) = \langle x - z, \vec{\rho}(x) - \vec{\rho}(z) \rangle = 0$$

implies $x = z$.

It follows that $\varphi(p, q) = x$ is a well-defined mapping from $\mathbb{R}^n_+ \times \mathbb{R}^n_+$ into $\mathbb{R}^n_+$. Suppose now that hypothesis $(H_k)$ holds, and therefore also each $\rho_i$ has positive derivative and is $\mathcal{C}^k$ for positive arguments. Let $W$ be an $n \times (m-1)$ matrix whose columns are a basis of $D$ (for instance, when $D = \mathcal{D}$, we may take the columns $b_1 - b_2, \ldots, b_1 - b_m$), and let $V$ be an $n \times (n - m + 1)$ matrix whose columns form a basis of $D^\perp$. Consider the map $F : \mathbb{R}^n_+ \to \mathbb{R}^n$ given by

$$F(x) := \begin{pmatrix} V'x \\ W'\vec{\rho}(x) \end{pmatrix}.$$

Observe that $F$ is of class $\mathcal{C}^k$. Denote by $J(x)$ the Jacobian of $F$ evaluated at an $x \in \mathbb{R}^n_+$.

We claim that $J(x)$ is nonsingular, for any $x$. The transpose of $J(x)$ can be expressed as a block matrix: $J(x)' = (V, TW)$, where $T = \mathrm{diag}(\rho'_1(\bar{x}_1), \ldots, \rho'_n(\bar{x}_n))$ is a symmetric positive definite matrix. As the columns of $TW$ are linearly independent, as are the columns of $V$, it will suffice to prove that the column spans of $V$ and of $TW$ intersect only at zero. Since the column spaces of $V$ and $W$ are orthogonal, this fact follows from the following observation: if $T$ is self-adjoint and positive definite, then $TW \bigcap W^\perp = \{0\}$ for any subspace $W$ of $\mathbb{R}^n$. (Proof: factor $T = R'R$, with $R$ nonsingular, and suppose that $Tw \in TW \bigcap W^\perp$. For any $w_0 \in W$, $0 = \langle Tw, w_0 \rangle = \langle Rw, Rw_0 \rangle$. Thus, $Rw \in RW \bigcap (RW)^\perp = \{0\}$, so $Rw = 0$, which implies $w = 0$.)

To conclude, we note that, given $p$ and $q$, $x = \varphi(p, q)$ is the unique solution of

$$G(x, p, q) = F(x) - \begin{pmatrix} V'p \\ W'\vec{\rho}(q) \end{pmatrix} = 0. \tag{32}$$

Thus, the Implicit Function Theorem gives that $\varphi$ is class $\mathcal{C}^k$, since $\frac{\partial G}{\partial x}(x, p, q) = J(x)$ has rank $n$ at all $x, p, q$ and $G$ is of class $\mathcal{C}^k$. ∎

**Remark 4.2** Note that, since $G(\varphi(p, q), p, q) \equiv 0$ in (32), we know that

$$J(\varphi(p, q)) \frac{\partial \varphi}{\partial p}(p, q) + \begin{pmatrix} V' \\ 0 \end{pmatrix} = 0 \tag{33}$$

for all $p, q$. Thus $(\partial \varphi / \partial p)(p, q) = -J(\varphi(p, q))^{-1}(V, 0)'$, and we conclude that $\partial \varphi / \partial p$ has constant rank, equal to $n - m + 1$. □

The following quantity measures deviations relative to $\mathcal{D}^\perp$. Let us define, for each $x, z \in \mathbb{R}^n_+$:

$$\delta(x, z) := \sum_{i=1}^m \sum_{j=1}^m \left( \langle b_i - b_j, \vec{\rho}(x) - \vec{\rho}(z) \rangle \right)^2. \tag{34}$$

**Remark 4.3** Note that $\delta(x, z) = 0$ if and only if $\vec{\rho}(x) - \vec{\rho}(z) \in \mathcal{D}^\perp$, since $\mathcal{D}$ is spanned by the differences $b_i - b_j$. Also, note that, since $x - x = 0 \in \mathcal{D}$, $x$ satisfies both (29) and (30) (which uniquely characterize $\varphi(x, z)$) if and only if $\vec{\rho}(x) - \vec{\rho}(z) \in \mathcal{D}^\perp$. To summarize:

$$x = \varphi(x, z) \Leftrightarrow \delta(x, z) = 0 \Leftrightarrow \vec{\rho}(x) - \vec{\rho}(z) \in \mathcal{D}^\perp \tag{35}$$

for any $x, z \in \mathbb{R}^n_+$. □

**Remark 4.4** We could also have defined a smaller, but basically equivalent, sum using only the generating differences $b_i - b_1$, $i = 1, \ldots, m-1$, but the above definition for $\delta$ seems more natural. Moreover, note that if we let

$$\underline{\delta}(x, z) := \sum_{i=1}^{m-1} (\langle b_i - b_1, \vec{\rho}(x) - \vec{\rho}(z) \rangle)^2 + \sum_{\ell=m}^{n} (\langle v_\ell, x - z \rangle)^2$$

where the $v_i$ constitute a basis of $\mathcal{D}^\perp$, then the uniqueness part of Lemma 4.1, applied with $D = \mathcal{D}$, gives that $\underline{\delta}(x, z) = 0$ if and only if $x = z$. (Because $x = \varphi(x, x)$, and $\underline{\delta}(x, z) = 0$ implies $z = \varphi(x, x)$.) □

**Remark 4.5** Recall the definition (23) of $\mathcal{R}$. We can extend $\varphi$ to a map $\mathcal{R} \times \mathbb{R}_+^n \to \mathbb{R}_+^n$ by defining $\varphi(p, q) := \varphi(p + d, q)$ for any $p \in \mathcal{R}$ and $q \in \mathbb{R}_+^n$, where $d \in \mathcal{D}$ is *any* element of $\mathcal{D}$ with the property that $p + d \in \mathbb{R}_+^n$. This definition is valid because, for any $d_1$ and $d_2$ in $\mathcal{D}$ such that $p + d_1$ and $p + d_2$ are in $\mathbb{R}_+^n$, we have $\varphi(p + d_1, q) = \varphi(p + d_2, q)$, since $x = \varphi(p + d_1, q)$ is uniquely determined by the two conditions $\delta(x, q) = 0$ and $x - (p + d_1) \in \mathcal{D}$, and the second condition is equivalent to $x - (p + d_2) \in \mathcal{D}$. Moreover, if hypothesis $(H_k)$ holds, then the extension of $\varphi$ to $\mathcal{R} \times \mathbb{R}_+^n$, which we denote again by $\varphi$, is of class $\mathcal{C}^k$. Indeed, given any $p_0 \in \mathcal{R}$, and a $d_0 \in \mathcal{D}$ so that $p_0 + d_0 \in \mathbb{R}_+^n$, it holds that $p + d_0 \in \mathbb{R}_+^n$ for all $p$ in some neighborhood of $p_0$, so the extension of $\varphi$ is obtained (on that neighborhood) by composition of the original $\varphi$ with the translation $(p, q) \mapsto (p + d_0, q)$, and the former is of class $\mathcal{C}^k$ by Lemma 4.1. □

**Remark 4.6** If hypothesis $(H_1)$ holds, then, for some continuous $d(\cdot) > 0$, and for all $x, z \in \mathbb{R}_+^n$,

$$x - z \in \mathcal{D} \Rightarrow \delta(x, z) \geq d(z) |x - z|^2 + o\left(|x - z|^2\right). \tag{36}$$

Indeed, $\delta(x, z) \geq \sum_{i=2}^{m} (\langle b_i - b_1, \vec{\rho}(x) - \vec{\rho}(z) \rangle)^2$, and $\vec{\rho}(x) = \vec{\rho}(z) + Q(x - z) + o(|x - z|)$, where $Q = Q(z) = \operatorname{diag}(\rho_1'(z_1), \ldots, \rho_n'(z_n))$ is the Jacobian of $\vec{\rho}$ evaluated at $z$ (nonsingular, by hypothesis $(H_1)$), so $\delta(x, z) \geq \sum_{i=2}^{m} (\langle b_i - b_1, x - z \rangle_Q)^2 + o(|x - z|^2)$, where $\langle \eta, \xi \rangle_Q := \langle \eta, Q\xi \rangle$ is a nondegenerate inner product on $\mathcal{D}$. This gives (36), because the $b_i - b_1$'s constitute a basis of $\mathcal{D}$. (Or, more concretely: $\delta(x, z) \geq \sum_{i=2}^{m} |W_z' \xi|^2 + o(|x - z|^2)$, where $\xi := R_z(x - z)$, $W_z := (R_z(b_2 - b_1), \ldots, R_z(b_m - b_1))$, and $R_z$ is a symmetric positive definite matrix, continuous on $z$, with $R_z^2 = Q(z)$. So, we may use that $|W_z' \xi| \geq d_0(z) |\xi|$, where $d_0(z)$ is the smallest singular value of $W_z'$ restricted to $R_z \mathcal{D}$, which is nonzero since $W_z' \xi = 0$ implies $\xi \in (R_z \mathcal{D})^\perp$, together with $|\xi| \geq d_1(z) |x - z|$, where $d_1(z)$ is the smallest eigenvalue of $R_z$, to obtain again an estimate (36).) □

## 5  Interior Equilibria

It is convenient to also express the dynamics (5) in matrix terms. Letting

$$\widetilde{A} = A - \operatorname{diag}\left(\sum_{i=1}^{m} a_{i1}, \ldots, \sum_{i=1}^{m} a_{im}\right) = \begin{pmatrix} -\sum_{i \neq 1}^{m} a_{i1} & a_{12} & \cdots & a_{1m} \\ a_{21} & -\sum_{i \neq 2}^{m} a_{i2} & \cdots & a_{2m} \\ \vdots & \vdots & \vdots & \vdots \\ a_{m1} & a_{m2} & \cdots & -\sum_{i \neq m}^{m} a_{im} \end{pmatrix}$$

we write:
$$\dot{x} = f(x) = B\widetilde{A}\Theta_B(x) \tag{37}$$

where $\Theta_B$ is the mapping

$$\Theta_B : \mathbb{R}^n \to \mathbb{R}_{\geq 0}^m : x \mapsto \left(e^{\langle b_1, \vec{\rho}(x) \rangle}, \ldots, e^{\langle b_m, \vec{\rho}(x) \rangle}\right)'$$

obtained as the composition of the maps $x \mapsto \vec{\rho}(x)$, $z \mapsto B'z$, and $y \mapsto (e^{y_1}, \ldots, e^{y_m})'$. In particular, since $\operatorname{rank} B = m$ and each $\rho_i$ maps $\mathbb{R}_+$ onto $\mathbb{R}$,

$$\text{the restriction } \Theta_B : \mathbb{R}_+^n \to \mathbb{R}_+^m \text{ is onto.} \tag{38}$$

Note that

$$\Theta_B(x) = \left(\theta_1(x_1)^{b_{11}} \theta_2(x_2)^{b_{21}} \ldots \theta_n(x_n)^{b_{n1}}, \ldots, \theta_1(x_1)^{b_{1m}} \theta_2(x_2)^{b_{2m}} \ldots \theta_n(x_n)^{b_{nm}}\right)'.$$

For any two $\bar{x}, x \in \mathbb{R}_+^n$, it holds that:

$$(\exists \kappa > 0) \; \Theta_B(x) = \kappa \Theta_B(\bar{x}) \iff \vec{\rho}(x) - \vec{\rho}(\bar{x}) \in \mathcal{D}^\perp \tag{39}$$

(recall the definition (10) of $\mathcal{D}$). To see this, denote $y := \Theta_B(x)$, $\bar{y} := \Theta_B(\bar{x})$. If $y = \kappa \bar{y}$, then, with $k := \ln \kappa$, $\ln y_j = k + \ln \bar{y}_j$, for each $j = 1, \ldots, m$. Thus $\langle b_j, \vec{\rho}(x) \rangle = k + \langle b_j, \vec{\rho}(\bar{x}) \rangle$, which implies that $\langle b_j, \vec{\rho}(x) - \vec{\rho}(\bar{x}) \rangle = k$ for all $j$, and therefore

$$\langle b_j - b_i, \vec{\rho}(x) - \vec{\rho}(\bar{x}) \rangle = 0 \quad \forall i, j.$$

Conversely, if this holds, we may define $k := \langle b_1, \vec{\rho}(x) - \vec{\rho}(\bar{x}) \rangle$, $\kappa := e^k$, and reverse all implications.

We also note, using once again that $B$ has full column rank, that $f(\bar{x}) = 0$ is equivalent to $\widetilde{A} \Theta_B(\bar{x}) = 0$, that is:

**Lemma 5.1** A state $\bar{x}$ is an equilibrium if and only if $\Theta_B(\bar{x}) \in \ker \widetilde{A}$. $\square$

We now consider the matrix $\widetilde{A}$. The row vector $\underline{1} = (1, \ldots, 1)$ has the property that $\underline{1} \widetilde{A} = (0, \ldots, 0)$, so, in particular, $\widetilde{A}$ is singular. The following is a routine consequence of the Perron-Frobenius (or finite dimensional Krein-Rutman) Theorem.

**Lemma 5.2** There exists $\bar{y} \in \mathbb{R}_+^m \cap \ker \widetilde{A}$ so that $(\mathbb{R}_{\geq 0}^m \setminus \{0\}) \cap \ker \widetilde{A} = \{\kappa \bar{y}, \; \kappa > 0\}$.

*Proof.* If $y \in \mathbb{R}_{\geq 0}^m$ is any eigenvector of $\widetilde{A}$, corresponding to an eigenvalue $\lambda$, it follows that $0 = \underline{1} \widetilde{A} y = \underline{1} \lambda y = \lambda q$, where $q := \underline{1} y$ is a positive number (because $y$, being an eigenvector, is nonzero), and therefore necessarily $\lambda = 0$. In other words, a nonnegative eigenvector can only be associated to the zero eigenvalue. Pick now any $\gamma > 0$ large enough such that all entries of $\widehat{A} := \widetilde{A} + \gamma I$ are nonnegative. Since the incidence graph $G(\widehat{A})$ coincides with $G(A)$, it follows that $\widehat{A}$ is also irreducible. By the Perron-Frobenius Theorem, the spectral radius $\sigma$ of $\widehat{A}$ is positive and it is an eigenvalue of $\widehat{A}$ of algebraic multiplicity one, with an associated positive eigenvector $\bar{y} \in \mathbb{R}_+^m$. Moreover, every nonnegative eigenvector $y \in \mathbb{R}_{\geq 0}^m$ associated to $\sigma$ is a positive multiple of $\bar{y}$. As adding $\gamma I$ moves eigenvalues by $\gamma$ while preserving eigenvectors (that is, $(\widetilde{A} + \gamma I) y = (\lambda + \gamma) y$ is the same as $\widetilde{A} y = \lambda y$), $\bar{y}$ is a positive eigenvector of the original matrix $\widetilde{A}$. It is necessarily in the kernel of $\widetilde{A}$, since we already remarked that any nonnegative eigenvector must be associated to zero. Finally, if $y$ is any other nonnegative eigenvector of $\widetilde{A}$, and in particular any element of $(\mathbb{R}_{\geq 0}^m \setminus \{0\}) \cap \ker \widetilde{A}$, then it is also a nonnegative eigenvector of $\widehat{A}$, and thus it must be a positive multiple of $\bar{y}$, completing the proof. ∎

**Corollary 5.3** The set of positive equilibria $E_+$ is nonempty. Moreover, pick any fixed $\bar{x} \in E_+$. Then, for any positive vector $x \in \mathbb{R}_+^n$, the following equivalence holds:

$$x \in E_+ \iff \delta(x, \bar{x}) = 0. \tag{40}$$

*Proof.* By Lemma 5.2, there is some $\bar{y} \in \mathbb{R}_+^m$ in $\ker \widetilde{A}$. By (38), there is some $\bar{x} \in \mathbb{R}_+^n$ such that $\Theta_B(\bar{x}) = \bar{y}$. In view of Lemma 5.1, $\bar{x}$ is an equilibrium.

Now fix any $\bar{x} \in E_+$ and any $x \in \mathbb{R}_+^n$, and let $y := \Theta_B(x)$, $\bar{y} := \Theta_B(\bar{x})$. Suppose $x \in E_+$. By Lemma 5.1, $y \in \ker \widetilde{A}$. By Lemma 5.2, every two positive eigenvectors of $\widetilde{A}$ are multiples of each other, so there is some $\kappa \in \mathbb{R}_+$ such that $y = \kappa \bar{y}$. By (39), $\vec{\rho}(x) - \vec{\rho}(\bar{x}) \in \mathcal{D}^\perp$. Conversely, if this holds, then, again by (39), $y = \kappa \bar{y}$. hence $y$ is also an eigenvector of $\widetilde{A}$, so by Lemma 5.1 we conclude $x \in E_+$. ∎

**Remark 5.4** Even if the assumption that each $\theta_i|_{\mathbb{R}_{\geq 0}}$ is onto $\mathbb{R}_{\geq 0}$ is dropped, the second part of Corollary 5.3 is still valid, as the proof of the previous lemmas did not use ontoness. We are not assured that equilibria exist, but if there is any equilibrium $\bar{x}$ then (40) holds. $\square$

**Remark 5.5** It is possible to generalize many of the results to the case when $A$ is not irreducible but, instead, there exists a permutation matrix $P$ with the property that $P'AP$ is a block matrix with irreducible blocks. Let us sketch this next. Assuming already such a reordering has taken place, the family of systems considered is as follows. The dynamics are described by

$$\dot{x} = f_1(x) + \ldots + f_L(x)$$

where $L$ is a positive integer (the "number of linkage classes"), and each $f_s$ has the form in Equation (5), for some matrices $A = A_s$ and $B = B_s$. We suppose that each $A_s$ is an irreducible nonnegative matrix of size $m_s \times m_s$, and each $B_s$ has size $n \times m_s$. We let $m = m_1 + \ldots + m_s$. Further, we also assume that the matrix $B$ which is composed from the blocks $B_s$, namely $B = (B_1, \ldots, B_L)$, satisfies the properties (7), (8), (9) required of $B$. We let $\mathcal{D}_s$ be the subspace of $\mathbb{R}^n$ spanned by the differences $b_i - b_j$ of columns of $B_s$, for each $s = 1, \ldots, L$, and now define $\mathcal{D}$ as the sum of the spaces $\mathcal{D}_s$. (Note that the column spaces of the $B_s$'s intersect only at zero, because we are assuming that the columns of $B$ are linearly independent; thus $\mathcal{D}$ is also a direct sum of the $\mathcal{D}_s$'s. The dimension of each $\mathcal{D}_s$ is $m_s - 1$, and $\mathcal{D}$ has dimension $m - L$.) Since the differences $b_i - b_j$ of columns of different $B_s$'s do not ever appear in the vector fields defining the system, the same argument as before shows that the cosets $p + \mathcal{D}$ (with the new definition of $\mathcal{D}$) are invariant. Classes are now defined using this $\mathcal{D}$, and are also invariant. We can express the dynamics in the form (37), where $\widetilde{A}$ is formed as before, starting from the $m \times m$ matrix $A$ that is obtained by using $A_1, \ldots, A_L$ as diagonal blocks.

Since each $\widetilde{A}_s$ is irreducible, we can find positive eigenvectors for the entire matrix $\widetilde{A}$, and the possible such eigenvectors are of the form $(\kappa_1 \bar{y}'_1, \ldots, \kappa_L \bar{y}_L)'$, where each $\bar{y}_s$ is a positive eigenvector of $A_s$ and the $\kappa_s$'s are positive numbers. The function $\Theta_B$ is still onto, so we obtain the existence of positive equilibria. Note that $\Theta_B(x)$ decomposes into blocks $\Theta_B^1, \ldots, \Theta_B^L$, mapping into $\mathbb{R}_{\geq 0}^{m_s}$ respectively, and (39) generalizes to: $\Theta_B^s(x) = \kappa \Theta_B^s(\bar{x})$ for some $\kappa > 0$ iff $\vec{\rho}(x) - \vec{\rho}(\bar{x}) \in \mathcal{D}_s^\perp$. Thus, uniqueness in each class holds just as before, since $\mathcal{D}^\perp$ is the intersection of the $\mathcal{D}_s^\perp$, $s = 1, \ldots, L$, and also $\bar{x}$ is an equilibrium iff all $f_i(\bar{x}) = 0$.

As a last comment along these lines, we remark that the assumption that $B$ has full rank can also be slightly relaxed, as follows. Suppose that the column spaces of the $B_i$'s intersect only at the origin, and the column space of some $B_s$ spans an affine space of dimension $m_s - 1$, i.e., $\mathcal{D}_s$ has dimension $m_s - 1$, but $B_s$ itself has rank $m_s - 1$ (instead of $m_s$). Then, we may add a state variable to the system, which satisfies a differential equation $\dot{x}_0 = 0$, and extend the system in such a manner that $B_s$ has rank $m_s$: this is equivalent to adding a row constantly equal to 1 to the matrix $B_s$ (and to the remaining matrices as well). The original system appears as the restriction of the new system to the invariant subset consisting of those states whose last coordinate is equal to one. $\square$

## 5.1 Proof of Part a in Theorem 2 and Theorem 6

Fix an $\bar{x} \in E_+$ and any $p \in \mathbb{R}_+^n$. Applied with $q = \bar{x}$ (and $D = \mathcal{D}$), Lemma 4.1 says that $x = \varphi(p, \bar{x})$ is uniquely characterized by (a) $x \in p + \mathcal{D}$ and (b) $\delta(x, \bar{x}) = 0$, the second of which is in turn equivalent by (40) to $x \in E_+$. In conclusion, $x = \varphi(p, \bar{x})$ is the unique point in $\mathbb{R}_+^n$ which is an equilibrium and lies in the same class as $p$. This proves Part a in Theorem 2, and, moreover, shows that

$$\pi(p) = \varphi(p, \bar{x}) \tag{41}$$

is the assignment of $p$ into the unique positive equilibrium in its class. Furthermore, by Lemma 4.1, $\varphi(\cdot, \bar{x})$ is of class $\mathcal{C}^k$ provided hypothesis ($H_k$) holds, which proves Theorem 6. ∎

**Remark 5.6** Theorem 2 is valid even under the weaker assumption that condition (28) holds, not necessarily assuming that each $\theta_i$ is onto $\mathbb{R}_{\geq 0}$. More precisely, *if* there is at least one equilibrium, then the conclusions hold. The stability parts will be remarked upon later, but regarding Part a, it is only necessary to note that Lemma 4.1 was proved under the weaker hypothesis, and that (40) holds (cf. Remark 5.4).

Now let us introduce the following condition, expressed in terms of the transpose of the vector $\underline{1} = (1, \ldots, 1)$ and the image of the transpose of $B$:

$$\underline{1}' \in B'(\mathbb{R}_+^n). \tag{42}$$

(Observe that, if the system is "homogeneous" in the sense that the exponents in each term add to a fixed constant, then (42) holds, since homogeneity amounts to the requirement that $B'\underline{1}_n = \underline{1}$, where we denote by $\underline{1}_n$ the vector of length $n$ consisting of 1's.) We claim that if (42) holds (instead of unboundedness), then there exists at least one interior equilibrium. For this we may assume, without loss of generality, that $\sigma_i = 1$ for all $i$ (if this were not the case, we simply rescale all the functions $\theta_i$, adjusting the $A$ matrix accordingly; thus we rewrite the system equations as those of a system of the same form, for which $\sigma_i = 1$). It follows that the image of $\Theta_B$ is $\exp(B'(-\mathbb{R}_+^n))$, where "exp" means taking exponentials of each coordinate. Now take any $\bar{y} \in \mathbb{R}_+^m \cap \ker \widetilde{A}$ and let $v \in \mathbb{R}^m$ be so that $\bar{y} = \exp(v)$. As $B'(\mathbb{R}_+^n)$ is an open set, and using (42), there are $\lambda > 0$ and $q \in \mathbb{R}_+^n$ such that $\underline{1}' - (1/\lambda)v = B'q$, or equivalently $-\lambda \underline{1} + v = B'(-\lambda q)$ which means that $e^{-\lambda}\bar{y} = \Theta_B(x)$ for some $x \in \mathbb{R}_+^n$. As $e^{-\lambda}\bar{y} \in \ker \widetilde{A}$, $x$ is an equilibrium, as wanted. □

**Remark 5.7** Recall Remark 4.5. After fixing an $\bar{x} \in E_+$, we may extend $\pi$ to $\mathcal{R}$ by letting $\pi(p) := \varphi(p, \bar{x})$ for any $p \in \mathcal{R}$. As $\varphi(p, \bar{x}) := \varphi(p + d, \bar{x})$ for some $d$ so that $p + d \in \mathbb{R}_+^n$, and $\varphi(p + d, \bar{x}) \in E_+$ by (41), we have that $\pi(p) \in E_+$ for all $p \in \mathcal{R}$. Note also that $\pi(p) = p$ if $p \in E_+$, since in that case $p$ is itself the unique positive equilibrium in its class. Thus,

$$p = \pi(p) \iff p \in E_+ \tag{43}$$

for any $p \in \mathcal{R}$. Furthermore, as $\varphi(p, \bar{x}) - (p + d) = \varphi(p + d, \bar{x}) - (p + d) \in \mathcal{D}$ by definition of $\varphi$, we know that $\varphi(p, \bar{x}) \in p + \mathcal{D}$, so

$$p - \pi(p) \in \mathcal{D} \tag{44}$$

for all $p \in \mathcal{R}$. We also have that

$$\pi(p_1) = \pi(p_2) \iff p_1 - p_2 \in \mathcal{D} \tag{45}$$

for all $p_1, p_2 \in \mathcal{R}$, by (44) and because $\pi(p_2)$ is the unique positive equilibrium in $p_2 + \mathcal{D} = p_1 + \mathcal{D}$ when $p_1 - p_2 \in \mathcal{D}$. Finally, note that $\pi$ is of class $\mathcal{C}^k$ if hypothesis ($H_k$) holds, once more by Remark 4.5. □

## 5.2 Proof of Theorem 5

Suppose that hypothesis ($H_k$) holds. Then the restriction of each $\rho_i$ to $\mathbb{R}_+$ is of class $\mathcal{C}^k$, and so is its inverse, which is a map $\mathbb{R} \to \mathbb{R}^+$. With some abuse of notation, we will denote in this section the restriction of each $\rho_i$ to $\mathbb{R}_+$, and more generally the associated vector function $\vec{\rho}$, restricted to $\mathbb{R}_+^n$, with the same symbols. Note that $\vec{\rho}$ is a $\mathcal{C}^k$-diffeomorphism (class $\mathcal{C}^k$ and has an inverse which is also of class $\mathcal{C}^k$). Now we fix any equilibrium $\bar{x} \in E_+$, and in terms of it, define:

$$M : \mathbb{R}^n \to \mathbb{R}_+^n : y \mapsto \vec{\rho}^{-1}(y + \vec{\rho}(\bar{x})),$$

which is a $\mathcal{C}^k$-diffeomorphism since both $y \mapsto y + \vec{\rho}(\bar{x})$ and $\vec{\rho}^{-1}$ are $\mathcal{C}^k$-diffeomorphisms. As an illustration, in the standard setup we would have $\rho_i(r) = \ln r$ and $\rho_i^{-1}(s) = e^s$, so $M(y) = (\bar{x}_1 e^{y_1}, \ldots, \bar{x}_n e^{y_n})'$.

We claim that $M(\mathcal{D}^\perp) = E_+$. Once that this is shown, Theorem 5 will be established, since diffeomorphisms map submanifolds into submanifolds, and $\mathcal{D}^\perp$ has dimension $n - (m-1)$. Note that $x \in M(\mathcal{D}^\perp)$ if and only if $\vec{\rho}(x) = y + \vec{\rho}(\bar{x})$ for some $y \in \mathcal{D}^\perp$, that is, if and only if $\vec{\rho}(x) - \vec{\rho}(\bar{x}) \in \mathcal{D}^\perp$, so Corollary 5.3 gives the equivalence with $x \in E_+$. ■

## 5.3 Proof of Theorem 7

We will construct a preliminary transformation $\widehat{\Phi}$ of $\mathcal{R}$ into $\mathcal{D} \times E_+$. Since $\mathcal{D}$ is a linear space of dimension $m-1$, and since by Theorem 5 we know that $E_+$ is a manifold diffeomorphic to $\mathbb{R}^{n-m+1}$, the result will then easily follow. We assume that hypothesis (H$_k$) holds, and fix any $\bar{x} \in E_+$.

Recalling Remark 5.7, we define, for $x \in \mathcal{R}$:

$$\widehat{\Phi}(x) = (\widehat{\Phi}_1(x), \widehat{\Phi}_2(x)) := (x - \pi(x), \pi(x))$$

(see Figure 3). The function $\widehat{\Phi}$ is of class $\mathcal{C}^k$ as a map into $\mathbb{R}^n \times \mathbb{R}^n_+$, and, because of (43),

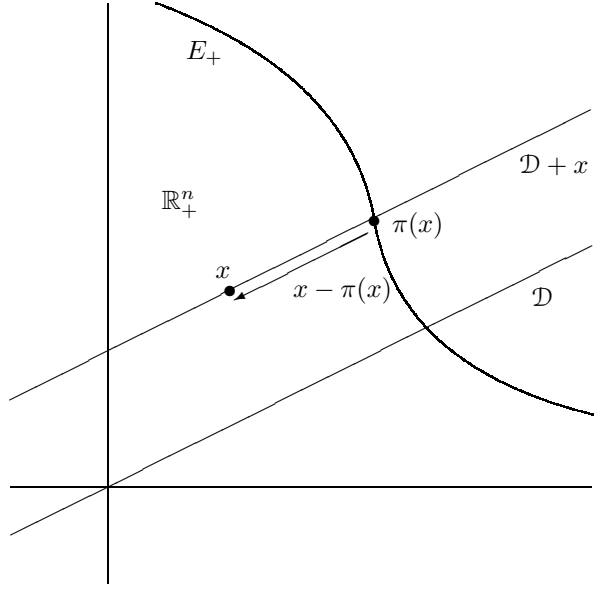

Figure 3: Change of coordinates

$$x \in E_+ \iff \widehat{\Phi}_1(x) = 0 \iff \widehat{\Phi}(x) = (0, x) \tag{46}$$

for all $x \in \mathcal{R}$. Also, property (44) gives that $\widehat{\Phi}_1(x) \in \mathcal{D}$ for all $x$, and we have that $\widehat{\Phi}_2(x) \in E_+$ (since $\pi$ assigns positive equilibria). So the image of $\widehat{\Phi}$ is contained in $\mathcal{D} \times E_+$, which is an embedded submanifold of $\mathbb{R}^n \times \mathbb{R}^n_+$ because $\mathcal{D}$ is a linear subspace and using Theorem 5. Thus, $\widehat{\Phi}$ defines also a $\mathcal{C}^k$ mapping into the submanifold $\mathcal{D} \times E_+$, and we view it as such.

Now consider the mapping

$$\Psi : \mathcal{D} \times \mathbb{R}^n_+ \to \mathbb{R}^n : (a, b) \mapsto a + b$$

and its image $\mathcal{R} = \Psi(\mathcal{D} \times \mathbb{R}^n_+) = \mathbb{R}^n_+ + \mathcal{D}$. If we show that the restriction $\widetilde{\Psi}$ of $\Psi$ to $\mathcal{D} \times E_+$ is the inverse of $\widehat{\Phi}$, then we will know that $\widehat{\Phi}$ is a diffeomorphism. It is clear from the formulas that $\widetilde{\Psi}(\widehat{\Phi}(x)) = x$. Pick now any $(a, b) \in \mathcal{D} \times E_+$, and let $x = \widetilde{\Psi}(a, b) = a + b$. Thus, $x - a = b \in E_+$, and $\varphi(b, \bar{x}) = b$ by (35) and (40). By the definition of $\varphi$, and thus $\pi$, on $\mathcal{R}$, we have that $\pi(x) = \varphi(x, \bar{x}) = \varphi(x - a, \bar{x})$, so $\widehat{\Phi}(x) = (x - \varphi(x - a, \bar{x}), \varphi(x - a, \bar{x}))$, which equals $(x - b, b) = (a, b)$. Thus also $\widehat{\Phi}(\widetilde{\Psi}(a, b)) = (a, b)$.

In conclusion, we have a diffeomorphism $\mathcal{R} \to \mathcal{D} \times E_+$, which can be composed with a diffeomorphism $\psi_1 \times \psi_2 : \mathcal{D} \times E_+ \to \mathbb{R}^{m-1} \times \mathbb{R}^{n-m+1}$ (where $\psi_1$ is a linear map) to yield the final diffeomorphism $\Phi = (\Phi_1, \Phi_2) : \mathcal{R} \to \mathbb{R}^{m-1} \times \mathbb{R}^{n-m+1}$. We let $F$ be the vector field obtained under this change of coordinates (namely, $F(X) = \Phi_*(\Phi^{-1}(X))f(\Phi^{-1}(X)))$ and write $F = (F_1, F_2)$ in block form. Note that, by (46), $\Phi_1(x) = 0$ if and only if $x \in E_+$. That is, block vectors $(X_1, X_2)$ with $X_1 = 0$ correspond under the coordinate change to $E_+$. That is, every vector of the form $(0, X_2)$ is an equilibrium. This implies that $F(0, X_2) = 0$ for all $X_2$, and in particular the same is true for the first block $F_1$. Also, by (45) we know that $x - z \in \mathcal{D}$ if and only if $\widehat{\Phi}_2(x) = \widehat{\Phi}_2(z)$, i.e., if and only if $\Phi_2(x) = \Phi_2(z)$.

Finally, we show that $F_2 \equiv 0$. It suffices to prove that $X_2(t) = \Phi_2(x(t))$ is constant, for an arbitrary solution $x(t)$ of (5) taking values in $\mathcal{R}$. This follows from the fact that classes are invariant, and hence the function $\varphi(x(t), \bar{x})$ is constant, which gives from the definition that $\widehat{\Phi}_2(x(t))$, and thus also $\Phi_2(x(t))$, is constant. ■

## 6 Boundary Equilibria

Fix any $x = (x_1, \ldots, x_n) \in \mathbb{R}^n_{\geq 0}$. We wish to study the implications of some coordinate $x_k$ vanishing. Recall the definition (13) of the sets $\mathcal{S}_j$. We will use repeatedly the following fact, for any $j \in \{1, \ldots, m\}$ and $k \in \{1, \ldots, n\}$:

$$x_k = 0 \text{ and } k \in \mathcal{S}_j \;\Rightarrow\; \theta_1(x_1)^{b_{1j}} \theta_2(x_2)^{b_{2j}} \ldots \theta_n(x_n)^{b_{nj}} = 0 \tag{47}$$

which is obvious, since $\theta_k(x_k)^{b_{kj}}$ vanishes when $x_k = 0$ and $b_{kj} \neq 0$. In particular,

$$x_k = 0 \;\Rightarrow\; \theta_1(x_1)^{b_{1j}} \theta_2(x_2)^{b_{2j}} \ldots \theta_n(x_n)^{b_{nj}} b_{kj} = 0 \tag{48}$$

since either $b_{kj} = 0$ or (47) applies. We state results for arbitrary perturbed systems (17). (Note that (49) below is merely the particular case of (50) that arises for systems with $g \equiv 0$, but it will be useful to have the conclusion for $f_k$ stated separately.)

**Lemma 6.1** Take any $i, j \in \{1, \ldots, m\}$ such that $a_{ij} \neq 0$. Then,

$$(\forall \ell \in \mathcal{S}_j)(x_\ell > 0) \;\Rightarrow\; (\forall k \in \mathcal{S}_i)(x_k > 0 \text{ or } f_k(x) > 0) \tag{49}$$

and

$$(\forall \ell \in \mathcal{S}_j)(x_\ell > 0) \;\Rightarrow\; (\forall k \in \mathcal{S}_i)(x_k > 0 \text{ or } f_k^*(x) > 0). \tag{50}$$

*Proof.* Pick any $k \in \mathcal{S}_i$, and assume that $x_k = 0$. The assumption "$\ell \in \mathcal{S}_j \Rightarrow x_\ell > 0$" is equivalent to $j \in J_{\ell,1} \Rightarrow x_\ell > 0$. So, since $x_k = 0$, necessarily $j \in J_{k,0}$. By (27), $f_k(x) = \beta_k(x)$, where $\beta_k$ is as in (26), and the index $j$ being considered does appear in the sum defining $\beta_k$. Moreover, for each $\ell \in \mathcal{S}_j$, $x_\ell \neq 0$ by hypothesis, so $\theta_1(x_1)^{b_{1j}} \theta_2(x_2)^{b_{2j}} \ldots \theta_n(x_n)^{b_{nj}} > 0$. On the other hand, $k \in \mathcal{S}_i$ means that $b_{ki} > 0$, and also $a_{ij} \neq 0$. Thus, the term involving this particular $i$ and $j$ in the sums defining $\beta_k(x)$ is positive (and the remaining terms are nonnegative). We conclude that $\beta_k(x) > 0$, establishing (49). In addition, since we are assuming that $x_k = 0$, property (27) also says that $f_k^*(x) = \beta_k(x) + g_k(x)$; as we proved that $\beta_k(x) > 0$, and $g_k(x) \geq 0$ (by assumption (18)), also (50) holds. ■

We will denote by $E^*$ (respectively, $E_+^*$ or $E_0^*$) the set of nonnegative (respectively, positive or boundary) equilibria of (17), i.e. the set of states $\bar{x} \in \mathbb{R}^n_{\geq 0}$ (respectively, $\in \mathbb{R}^n_+$ or $\in \mathbb{R}^n_0$), such that $f^*(\bar{x}) = 0$. When $g \equiv 0$, these are the same as the sets $E$, $E_+$, and $E_0$, and, more generally, for systems (17), we still use $E$, etc., for the equilibria of the system (5) which has the same $f$.

**Proposition 6.2** For any system (17), and for an arbitrary $x \in \mathbb{R}^n_{\geq 0}$, consider the following properties:

1. $x \in E_0^*$.

2. $(\forall j \in \{1, \ldots, m\}) (\exists k \in \mathcal{S}_j) \; x_k = 0$.

3. $(\forall j \in \{1, \ldots, m\}) \; \theta_1(x_1)^{b_{1j}} \theta_2(x_2)^{b_{2j}} \ldots \theta_n(x_n)^{b_{nj}} = 0$.

4. $x \in E_0$.

Then $1 \Rightarrow 2 \Leftrightarrow 3 \Rightarrow 4$.

*Proof.* [$1 \Rightarrow 2$] Pick any $x \in E_0^*$. If the second property is false, then there is some index $j$ such that $x_\ell > 0$ for all $\ell \in \mathcal{S}_j$. We claim that for every index $i$, $x_k > 0$ for all $k \in \mathcal{S}_i$. Since $\bigcup_j \mathcal{S}_j = \{1, \ldots, n\}$ (recall hypothesis (9)), this will mean that $x_k > 0$ for all $k$, so $x$ could not have been a boundary point, a contradiction. Let $J = \{j \,|\, x_\ell > 0 \,\forall \ell \in \mathcal{S}_j\}$ and let $I = \{1, \ldots, m\} \setminus J$. We know that $J \neq \emptyset$ and must prove that $I = \emptyset$. Suppose by contradiction that $I \neq \emptyset$. Pick some $i \in I$ and $j \in J$ such that $a_{ij} \neq 0$ (irreducibility of $A$), and take any $k \in \mathcal{S}_i$. We claim that $x_k > 0$. Suppose that this is not the case, i.e. $x_k = 0$. From (50) in Lemma 6.1, we conclude that $f_k^*(x) > 0$. This contradicts the fact that $x$ is in $E^*$. In conclusion, $x_k > 0$ for all $k \in \mathcal{S}_i$, contradicting the fact that $i \in I$.

[$2 \Leftrightarrow 3$] The product in 3 can vanish if and only if some $b_{kj} > 0$ and $x_k = 0$ for the same $k$, but this is precisely the condition in 2.

[$3 \Rightarrow 4$] Since all terms in the definition (5) of $f(x)$ vanish, it follows that $x \in E$. On the other hand, $x$ must be a boundary point, since otherwise no product as in 3 could vanish. ∎

For systems (5), all the properties are equivalent, since $E_0^* = E_0$. More generally, this stronger result holds for all systems of form (21):

**Proposition 6.3** *For any system (21) and for an arbitrary $x \in \mathbb{R}_{\geq 0}^n$, the properties in Proposition 6.2 are equivalent.*

*Proof.* We prove that $4 \Rightarrow 2 \Rightarrow 1$. If 4 holds, we apply Proposition 6.2 to the system (5) which has the same $f$, to conclude that 2 holds. Assume now that 2 (and hence also 4) holds. Since we already know that $f(x) = 0$, all that needs to be verified is that $\Delta_{ij}(x) = 0$ for all $i, j$. Pick any $i, j$. Choose $k \in \mathcal{S}_j$ so that $x_k = 0$. Property (20) then gives the conclusion. ∎

## 7  Invariance

We start with an easy but key observation.

**Lemma 7.1** *Suppose that $x : [0, t^*] \to \mathbb{R}_{\geq 0}^n$ is any solution of (17). Then the following implication holds for any $k = 1, \ldots, n$:*
$$x_k(0) > 0 \quad \Rightarrow \quad x_k(t^*) > 0\,.$$

*Proof.* Suppose that $k$ is so that $x_k(0) > 0$. Let $F : \mathbb{R}^2 \to \mathbb{R}$ be the (locally Lipschitz) function which coincides for $t \in [0, t^*]$ and $y \in \mathbb{R}$ with
$$F(t, y) := f_k^*(x_1(t), \ldots, x_{k-1}(t), y, x_{k+1}(t), \ldots, x_n(t))$$
and has $F(t, y) = F(0, y)$ for $t < 0$ and $F(t, y) = F(t^*, y)$ for $t > t^*$. Note that $F(t, 0) \geq 0$ for all $t$, because (27) says that $f_k^*(x) = \widetilde{\beta}_k(x) \geq 0$ when $x_k = 0$. For $t \in [0, t^*]$, the scalar function $y(t) := x_k(t)$ satisfies the scalar differential equation
$$\dot{y}(t) \;=\; F(t, y(t))\,.$$

We must prove that $y$ never vanishes. For this, we let $G(t, p) := F(t, p) - F(t, 0)$ and introduce the auxiliary initial value problem
$$\dot{z}(t) \;=\; G(t, z(t)), \quad z(0) = y(0)\,.$$

Since $G$ is locally Lipschitz (in both variables, but just on $z$ locally uniformly on $t$ would be enough) and 0 is an equilibrium of $\dot{z} = G(t, z)$, $z(t) > 0$ for all $t$ in its domain of definition. Moreover, we have that $\dot{z} = G(t, z) \leq F(t, z)$ for all $t$. By a standard comparison theorem, see e.g. [15], Corollary I.6.3, we know that $z(t) \leq y(t)$ for all $t$ in the common domain of definition of $z$ and $y$. Since $y(t^*)$ is well-defined, $z(t)$ remains bounded, and thus is defined as well for $t = t^*$. So, $y(t^*) \geq z(t^*) > 0$. ∎

**Corollary 7.2** The set $\mathbb{R}^n_+$ is forward invariant for (17).

*Proof.* Consider any solution $x : [0, T] \to \mathbb{R}^n$ of (17), and suppose that $x(0) \in \mathbb{R}^n_+$. We must prove that $x(T) \in \mathbb{R}^n_+$. Since $x(0)$ is in the interior of $\mathbb{R}^n_{\geq 0}$, the only way that the conclusion could fail is if $x(t) \in \mathbb{R}^n_0$ for some $t \in (0, T]$. We assume that this happens and derive a contradiction. Let $t^* := \min\{t \in [0, T] \,|\, x(t) \in \mathbb{R}^n_0\} > 0$. By minimality, for all $i$, $x_i(t) > 0$ for all $t \in [0, t^*)$, and in particular $x_i(t) \in \mathbb{R}^n_{\geq 0}$ for all $t \in [0, t^*]$, and also there is some index $k$ such that $x_k(t^*) = 0$. But this contradicts the conclusion of Lemma 7.1. ∎

The closure of an invariant set is also invariant, so:

**Corollary 7.3** The set $\mathbb{R}^n_{\geq 0}$ is forward invariant for (17).

Note that Corollaries 7.2 and 7.3 prove Lemma 2.1 as well as the invariance parts of Theorems 3 and 4.

**Lemma 7.4** Consider any solution $x : [0, T] \to \mathbb{R}^n$ of (17) for which $x(0) \in \mathbb{R}^n_{\geq 0}$. Suppose that there is some $j_0 \in \{1, \ldots, m\}$ such that

$$(\forall \ell \in \mathcal{S}_{j_0}) \; [\, x_\ell(0) > 0 \;\; \text{or} \;\; f^*_\ell(x(0)) > 0 \,].$$

Then, $x(t) \in \mathbb{R}^n_+$ for all $t \in (0, T]$.

*Proof.* We know that $x(t) \in \mathbb{R}^n_{\geq 0}$, by Corollary 7.3, and we need to prove that $x_k(t) > 0$ for all $k$ and all $t \in (0, T]$. Let

$$I := \{i \in \{1, \ldots, m\} \,|\, x_k(t) > 0 \; \forall k \in \mathcal{S}_i, \forall t \in (0, T]\}.$$

Since $\bigcup_i \mathcal{S}_i = \{1, \ldots, n\}$, it will be enough to show that $I = \{1, \ldots, m\}$.

We start by remarking that $j_0 \in I$, because, for each $k \in \mathcal{S}_{j_0}$, either $x_k(0) > 0$, and then Lemma 7.1, applied on any subinterval $[0, t^*]$, says that $x_k(t) > 0$ for all $t$, or $\dot{x}_k(0) = f^*_k(x(0)) > 0$ and $x_k(0) = 0$ imply that $x_k(t) > 0$ for all $t$ small enough, so also (again by Lemma 7.1) for all $t$.

Suppose by way of contradiction that $H := \{1, \ldots, m\} \setminus I \neq 0$. Pick some $i \in I$ and $h \in H$ so that $a_{hi} \neq 0$ (irreducibility of $A$). We will show that, for any given $t_0 \in (0, T]$, and for any given $k \in \mathcal{S}_h$, $x_k(t_0) > 0$, and this will contradict $h \in H$.

Since $i \in I$, $x_\ell(t_0/2) > 0$ for all $\ell \in \mathcal{S}_i$. Then, we can apply Lemma 6.1, to obtain that $x_k(t_0/2) > 0$ or $f^*_k(x(t_0/2)) > 0$. As before, if $x_k(t_0/2) > 0$ then via Lemma 7.1 we conclude the positivity of $x_k(t)$ for all $t > t_0/2$, and in particular of $x_k(t_0)$. And if instead $f_k(x(t_0/2)) > 0$, then $\dot{x}_k(t_0/2) > 0$ and $x_k(t_0/2) \geq 0$ imply $x_k(t) > 0$ for all $t > t_0/2$ near $t_0/2$, and so by Lemma 7.1 once again $x_k(t_0) > 0$. ∎

For the special case of systems (5), and more generally perturbed systems (21), we conclude that every trajectory starting on the boundary $\mathbb{R}^n_0$ which is not an equilibrium must immediately enter the positive orthant:

**Corollary 7.5** Consider any solution $x : [0, T] \to \mathbb{R}^n_{\geq 0}$ of (21) for which $x(0) \notin E^*_0$. Then, $x(t) \in \mathbb{R}^n_+$ for all $t \in (0, T]$.

*Proof.* By Proposition 6.3, $x(0) \notin E_0^*$ implies that there is some $j \in \{1, \ldots, m\}$ such that $x_\ell(0) \neq 0$ for all $\ell \in \mathcal{S}_j$. So, Lemma 7.4 insures that $x(t) \in \mathbb{R}_+^n$ for all $t \in (0, T]$. ∎

Another special case of that of the feedback system (16), and more generally when property (19) holds. This property says that there is some $j_0 \in \{1, \ldots, m\}$ such that, for every $\ell \in \mathcal{S}_{j_0}$, it holds for each $z \in \mathbb{R}_{\geq 0}^n$ that either (a) $z_\ell > 0$, or (b) $z_\ell = 0$ and $g_\ell(z) > 0$ (in which case also $f_k^*(z) = \beta_k(z) + g_k(z) > 0$). Applying with $z = x(0)$ and using Lemma 7.4 with this $j_0$, we have:

**Corollary 7.6** Suppose that property (19) holds. and consider any solution $x : [0, T] \to \mathbb{R}_{\geq 0}^n$ of (17). Then, $x(t) \in \mathbb{R}_+^n$ for all $t \in (0, T]$. ∎

We also note, for further reference:

**Corollary 7.7** Consider a system (21), and pick any $\xi \in \mathbb{R}_{\geq 0}^n$. Then, either $\xi \in E_0^*$ or $\xi$ belongs to some positive class.

*Proof.* If $\xi \notin E_0^*$, then Corollary 7.5, applied to a solution starting from $\xi$, and forward invariance of classes, give that $\xi$ is in the positive class containing $x(T)$. ∎

## 8 Stability

We start by establishing some useful estimates.

**Lemma 8.1** Define the following quadratic function:

$$Q(\eta_1, \ldots, \eta_m) := \sum_{i=1}^{m} \sum_{j=1}^{m} a_{ij}(\eta_i - \eta_j)^2. \tag{51}$$

Then, there exists a constant $\kappa > 0$ such that

$$Q(q_1, \ldots, q_m) \geq \kappa \sum_{i=1}^{m} \sum_{j=1}^{m} (q_i - q_j)^2 \tag{52}$$

for all $(q_1, \ldots, q_m) \in \mathbb{R}^m$.

*Proof.* We first observe that

$$Q(q_1, \ldots, q_m) = 0 \Rightarrow q_i = q_m, i = 1, \ldots, m - 1.$$

Indeed, obviously $Q(q_1, \ldots, q_m) = 0$ implies $q_i = q_j$ for each pair $i, j$ for which $a_{ij} \neq 0$. Now let $I$ be the set of indices $i$ such that $q_i = q_m$, and $J$ its complement; as $m \in I$, $I \neq \emptyset$. We need to see that $J = \emptyset$. Suppose that $J \neq \emptyset$. The connectedness of the incidence graph of $A$ provides an $i \in I$ and $j \in J$ such that $a_{ij} \neq 0$. Thus, $q_j = q_i = q_m$, contradicting $j \in J$.

Let us introduce next a quadratic form in $m - 1$ variables:

$$\sum_{i=1}^{m-1}\sum_{j=1}^{m-1} a_{ij}(\xi_i - \xi_j)^2 + \sum_{i=1}^{m-1} a_{im}\xi_i^2 + \sum_{j=1}^{m-1} a_{mj}\xi_j^2,$$

which we denote as $P(\xi_1, \ldots, \xi_{m-1})$. Since $(\eta_i - \eta_m) - (\eta_j - \eta_m) = \eta_i - \eta_j$ for all $i, j$, one has

$$Q(\eta_1, \ldots, \eta_m) = P(\eta_1 - \eta_m, \ldots, \eta_{m-1} - \eta_m).$$

Note that $P$ is positive definite: if $P(q_1, \ldots, q_{m-1}) = 0$, then $Q(q_1, \ldots, q_{m-1}, 0) = 0$, which as already observed implies that all $q_i = 0$. Thus, there is some constant $\kappa_0 > 0$ such that

$$P(p_1, \ldots, p_{m-1}) \geq \kappa_0 \sum_{i=1}^{m-1} p_i^2$$

for all $(p_1, \ldots, p_{m-1}) \in \mathbb{R}^{m-1}$, which means that

$$Q(q_1, \ldots, q_m) \geq \kappa_0 \sum_{i=1}^{m-1} (q_i - q_m)^2 \tag{53}$$

for all $(q_1, \ldots, q_m) \in \mathbb{R}^m$. As $(q_i - q_j)^2 \leq 2(q_i - q_m)^2 + 2(q_j - q_m)^2$ for all $i, j$, we may re-express the estimate (53) in the form (52), using a smaller constant $\kappa$ which depends only on $\kappa_0$ and $m$. ∎

The following estimate will be the basis of a Lyapunov function property to be established later.

**Lemma 8.2** There exist two continuous functions

$$v : \mathbb{R}^n \to \mathbb{R}^n, \quad c : \mathbb{R}_+^n \to \mathbb{R}_+$$

such that, for every pair of points $x, z$ in $\mathbb{R}_+^n$:

$$\langle \vec{\rho}(x) - \vec{\rho}(z), f(x) \rangle \leq -c(z) \delta(x, z) + \langle v(\vec{\rho}(x) - \vec{\rho}(z)), f(z) \rangle. \tag{54}$$

*Proof.* As $B$ has full column rank, there is an $m \times n$ matrix $B^\#$ (for instance, its pseudo-inverse) such that $B^\# B = I$. We let

$$v(\sigma) := \left( (e^{\langle b_1, \sigma \rangle}, \ldots, e^{\langle b_m, \sigma \rangle}) B^\# \right)'$$

and

$$c_0(\zeta) := \min_{j=1,\ldots,m} e^{\langle b_j, \vec{\rho}(\zeta) \rangle}.$$

Now take any pair of positive vectors $x, z$. Denote, for each $j = 1, \ldots, m$:

$$q_j := \langle b_j, \vec{\rho}(x) - \vec{\rho}(z) \rangle$$

and observe that

$$\langle b_i, v(\vec{\rho}(x) - \vec{\rho}(z)) \rangle = e^{q_i}, \quad i = 1, \ldots, m$$

so, using formula (24),

$$\langle v(\vec{\rho}(x) - \vec{\rho}(z)), f(z) \rangle = \sum_{i=1}^{m} \sum_{j=1}^{m} a_{ij} e^{\langle b_j, \vec{\rho}(z) \rangle} (e^{q_i} - e^{q_j}). \tag{55}$$

Therefore (writing $g(x, z) = \langle v(\vec{\rho}(x) - \vec{\rho}(z)), f(z) \rangle$ for simplicity):

$$\langle \vec{\rho}(x) - \vec{\rho}(z), f(x) \rangle = \sum_{i=1}^{m} \sum_{j=1}^{m} a_{ij} e^{\langle b_j, \vec{\rho}(x) \rangle} (q_i - q_j)$$

$$= \sum_{i=1}^{m} \sum_{j=1}^{m} a_{ij} e^{\langle b_j, \vec{\rho}(z) \rangle} e^{q_j} (q_i - q_j)$$

$$= \sum_{i=1}^{m} \sum_{j=1}^{m} a_{ij} e^{\langle b_j, \vec{\rho}(z) \rangle} (e^{q_j} (q_i - q_j) - (e^{q_i} - e^{q_j})) + g(x, z) \tag{56}$$

$$\leq -\frac{1}{2} \sum_{i=1}^{m} \sum_{j=1}^{m} a_{ij} e^{\langle b_j, \vec{\rho}(z) \rangle} (q_i - q_j)^2 + g(x, z) \tag{57}$$

$$\leq -\frac{1}{2} c_0(z) \sum_{i=1}^{m} \sum_{j=1}^{m} a_{ij} (q_i - q_j)^2 + g(x, z)$$

$$= -\frac{1}{2} c_0(z) Q(q_1, \ldots, q_m) + g(x, z),$$

where $Q$ is the quadratic form in Lemma 8.1. Equality (56) follows by adding and subtracting $g(x, z)$ and using (55). To justify (57), we note first that, for each $a > 0$, the function $\mathbb{R}_{\geq 0} \to \mathbb{R}$:

$$f_a(r) := e^a(r - a) - e^r + e^a + \frac{1}{2}(r - a)^2$$

is always $\leq 0$ (because $f_a(0) = -e^a a - 1 + e^a + (1/2)a^2 < 0$, $f_a(r) \to -\infty$ as $r \to +\infty$, and $f'_a(r) = e^a - e^r + (r - a) \neq 0$ for all $r > 0$). Now we use the inequality $e^a(r - a) - e^r + e^a \leq -\frac{1}{2}(r - a)^2$ in each term of the sum with $a = q_j$ and $r = q_i$ (recall $a_{ij}e^{\langle b_j, \vec{\rho}(z)\rangle} \geq 0$).

Lemma 8.1 gives that $Q(q_1, \ldots, q_m) \geq \kappa \delta(x, z)$. Thus, we may take $c(z) := \kappa c_0(z)/2$. ∎

## 8.1 An Entropy Distance

We will show the stability conclusions even if one does not impose the assumption that each $\theta_i$ is onto $\mathbb{R}_{\geq 0}$. Note that none of the results in Sections 6 and 7 used this condition. Thus, we suppose that each $\theta_i : \mathbb{R} \to [0, \sigma_i)$ (where $0 < \sigma_i \leq \infty$) is locally Lipschitz, has $\theta_i(0) = 0$, satisfies $\int_0^1 |\rho_i(y)| \, dy < \infty$ (where $\rho_i(y) = \ln \theta_i(y)$), and its restriction to $\mathbb{R}_{\geq 0}$ is strictly increasing and onto $[0, \sigma_i)$, with $\sigma_i \leq +\infty$. We let $\bar{\rho}_i := \ln \sigma_i$. (For the stability results, we will not even need to ask for condition (28).)

For any fixed constant $c \in (-\infty, \bar{\rho}_i)$, and each $i = 1, \ldots, n$, we consider the following function:

$$R_i^c(r) := \int_1^r \rho_i(s) \, ds - cr \,.$$

This function is a well-defined continuous mapping $\mathbb{R}_{\geq 0} \to \mathbb{R}$, continuously differentiable for $r \in (0, \bar{\rho}_i)$. Moreover, $R_i^c$ achieves a global minimum at the unique $r^c = r_i^c \in \mathbb{R}_+$ where $\rho_i(r^c) = c$, decreases for $r \in [0, r^c]$, and increases to $+\infty$ for $r > r^c$.

The following function will play a central role:

$$W : \mathbb{R}_{\geq 0}^n \times \mathbb{R}_+^n \to \mathbb{R} : (x, z) \mapsto \sum_{i=1}^n R_i^{\rho_i(z_i)}(x_i) \,.$$

The above-mentioned properties of the functions $R_i^{\rho_i(z_i)}$ imply that

$$x \neq z \quad \Rightarrow \quad W(x, z) > W(z, z), \tag{58}$$

i.e., for each fixed $z \in \mathbb{R}_+^n$, the function $W(\cdot, z)$ has a unique global minimum, at $z$. Note also that the gradient of $W(\cdot, z)$:

$$\frac{\partial W}{\partial x}(x, z) = (\vec{\rho}(x) - \vec{\rho}(z))' \tag{59}$$

(defined for $x \in \mathbb{R}_+^n$) vanishes only at $x = z$ and that (since $R_i^{\rho_i(z_i)}(x_i) \to +\infty$ if $x_i \to +\infty$)

$$|x| \to +\infty \quad \Rightarrow \quad W(x, z) \to +\infty \tag{60}$$

for every given $z$. As $W(\cdot, z)$ is continuous, this implies that

$$\{x \mid W(x, z) \leq w\} \tag{61}$$

is compact for every $z$ and every $w \in \mathbb{R}$.

**Remark 8.3** In the standard setup $\rho_i = \ln$, $W(x, z) = \sum_{i=1}^n x_i \ln x_i - x_i - x_i \ln z_i$. Then this formula, when states $x$ are interpreted probabilistically in applications such as chemical networks, is suggested by "relative entropy" considerations. □

## 8.2 Main Stability Results

The next lemma will be applied later with $S$ equal to the whole space or to a class, depending on the type of system.

**Lemma 8.4** Let $S \subseteq \mathbb{R}^n_{\geq 0}$ be a closed set, and pick $\bar{x} \in S$. Suppose that

$$\langle \vec{\rho}(x) - \vec{\rho}(\bar{x}), f^*(x) \rangle < 0 \qquad (62)$$

is valid for all $x \in S \bigcap \mathbb{R}^n_+$, $x \neq \bar{x}$. Consider any $\xi \in \mathbb{R}^n_{\geq 0}$ for which the maximal solution $x(t)$ of (17) with $x(0) = \xi$ is included in $S$.

1. If the system has the form (21) then $x(t)$ is defined for all $t \geq 0$, and
$$x(t) \to E_0 \bigcup \{\bar{x}\} \text{ as } t \to +\infty\,.$$

2. If (19) holds, then $x(t)$ is defined for all $t \geq 0$, and
$$x(t) \to \bar{x} \text{ as } t \to +\infty\,.$$

Furthermore, in either case $\bar{x}$ is an equilibrium of (17), asymptotically stable relative to $S$.

*Proof.* Fix $S$, $\bar{x}$, and $\xi$ as in the statement, and let $x(\cdot)$ be the maximal solution of (17) with $x(0) = \xi$, defined a priori on some maximal interval $[0, t^*)$. We will use

$$V(x) := W(x, \bar{x}) - W(\bar{x}, \bar{x}) \qquad (63)$$

as a Lyapunov-like function. By (58), this function is positive definite relative to the equilibrium $\bar{x}$, i.e., $V(x) \geq 0$ for all $x \in \mathbb{R}^n_{\geq 0}$, and $V(x) = 0$ if and only if $x = \bar{x}$. Moreover, $V$ is proper, meaning that the sublevel sets $\{x \,|\, V(x) \leq w\}$ are compact, for all $w \in \mathbb{R}_{\geq 0}$, by (61). Finally, $V$ is continuously differentiable in the interior $\mathbb{R}^n_+$, and, using (59),

$$\nabla V(x)\, f^*(x) < 0 \qquad (64)$$

whenever $x \in S \bigcap \mathbb{R}^n_+$, $x \neq \bar{x}$, by (62).

If the system has the form (21) and $\xi \in E_0$, then $\xi \in E_0^*$ (cf. Proposition 6.3), and thus, in that special case, $x(t)$ is of course defined for all $t > 0$ and converges to an equilibrium in $E_0$ (it is constant, in fact).

So, from now on, we suppose that either the system has the form (21) and $\xi \notin E_0$, or that (19) holds. By, respectively, Corollary 7.5 and Corollary 7.6, we know that $x(t) \in S \bigcap \mathbb{R}^n_+$ for all $t \in [0, t^*)$. So $V(x(t))$ is differentiable for $t \in (0, t^*)$, and $dV(x(t))/dt \leq 0$ (by (64) if $x(t) \neq \bar{x}$, and obvious otherwise, since $\nabla V(\bar{x}) = 0$), which means that $V(x(t))$ is nonincreasing. Since $V$ is proper, this means that the maximal trajectory is precompact, and hence it is defined on the entire interval $[0, +\infty)$, as claimed. Furthermore, the LaSalle Invariance Theorem implies that

$$x(t) \to \Omega^+(\xi) \text{ as } t \to +\infty\,,$$

where $\Omega^+(\xi)$ is the $\omega$-limit set of $\xi$, which is a compact subset of $\{p \,|\, V(p) = a\}$, for some $a \geq 0$. As the set $S$ is closed, $\Omega^+(\xi) \subseteq S$. We pick any $\zeta \in \Omega^+(\xi)$, and show that necessarily $\zeta = \bar{x}$ or, in the case of systems (21), $\zeta \in E_0$.

If $\zeta = \bar{x}$, we are done, so we may assume from now on that $\zeta \neq \bar{x}$. Similarly, if the system is of type (21), we will assume that $\zeta \notin E_0$. We now derive a contradiction.

Consider the forward trajectory $z(t)$ starting from $\zeta$. Since $\Omega^+(\xi)$ is invariant and a subset of $S$, $z(t) \in S$ for all $t$, and $z(t) \neq \bar{x}$ for all $t$ (since otherwise $z(t) \equiv \bar{x}$). Moreover, $z(t) \in \mathbb{R}^n_+$ for all $t > 0$ (using either Corollary 7.5 or Corollary 7.6). Thus (64) says that $dV(z(t))/dt < 0$ for all $t > 0$, which means that $V(z(\cdot))$ is strictly decreasing. But this is a contradiction, since $V(z(t)) \equiv a$.

The stability statement is a simple consequence of the fact that $V$ is a Lyapunov function (see e.g. [21], section 5.7) relative to $\bar{x}$ for the dynamics restricted to $S \bigcap \mathbb{R}^n_+$. ∎

### 8.2.1 Proof of Theorem 3

Suppose $r = n - m + 1$, $\gamma_1, \ldots, \gamma_r$ are positive, $k_1, \ldots, k_r \in \{1, \ldots, n\}$ are such that (14) holds, and $j$ is such that (15) holds. Fix any equilibrium $\bar{x} \in E_+$.

We must prove that, for any given $\xi \in \mathbb{R}^n_{\geq 0}$, the maximal solution $x(t)$ of (16) with $x(0) = \xi$ is defined for all $t \geq 0$ (the fact that the solution remains in $\mathbb{R}^n_{\geq 0}$ has already been established) and that $\bar{x}$ is a globally asymptotically stable equilibrium of (16).

We apply Lemma 8.4, with $S = \mathbb{R}^n_{\geq 0}$. Notice that

$$\langle \vec{\rho}(x) - \vec{\rho}(\bar{x}), f(x) \rangle \leq -c(\bar{x})\, \delta(x, \bar{x}) \leq 0 \tag{65}$$

for all $x \in \mathbb{R}^n_+$ by (54) (since $f(\bar{x}) = 0$). The inequality is strict unless $\delta(x, \bar{x}) = 0$, i.e.,

$$\vec{\rho}(x) - \vec{\rho}(\bar{x}) \in \mathcal{D}^\perp. \tag{66}$$

On the other hand,

$$\begin{aligned}
\langle \vec{\rho}(x) - \vec{\rho}(\bar{x}), g(x) \rangle &= \sum_{\ell=1}^r \gamma_\ell (\bar{x}_{k_\ell} - x_{k_\ell}) \langle \vec{\rho}(x) - \vec{\rho}(\bar{x}), e_{k_\ell} \rangle \\
&= \sum_{\ell=1}^r \gamma_\ell (\bar{x}_{k_\ell} - x_{k_\ell}) (\rho_{k_\ell}(x_{k_\ell}) - \rho_{k_\ell}(\bar{x}_{k_\ell})) \leq 0
\end{aligned}$$

where the last inequality follows from the fact that each $\rho_i$ is an increasing function, and the inequality is strict unless

$$\langle \vec{\rho}(x) - \vec{\rho}(\bar{x}), e_{k_\ell} \rangle = \rho_{k_\ell}(x_{k_\ell}) - \rho_{k_\ell}(\bar{x}_{k_\ell}) = 0 \tag{67}$$

for $\ell = 1 \ldots r$. Thus $\langle \vec{\rho}(x) - \vec{\rho}(\bar{x}), f^*(x) \rangle$ equals

$$\langle \vec{\rho}(x) - \vec{\rho}(\bar{x}), f(x) \rangle + \langle \vec{\rho}(x) - \vec{\rho}(\bar{x}), g(x) \rangle \leq 0, \tag{68}$$

and this inner product can only vanish if both (66) and (67) hold, which implies, because of (14), that $\vec{\rho}(x) = \vec{\rho}(\bar{x})$, i.e., $x = \bar{x}$. This means that (62) holds. As the system (16) is a system (17) for which (19) holds, Lemma 8.4 provides the proof of global stability. ∎

### 8.2.2 Proof of Theorem 1

Consider any system (21) and any maximal solution $x(\cdot)$. By Corollary 7.7, either $x(t) \equiv \xi \in E_0^* = E_0$, in which case of course $x(t) \to E$, or there is a positive class $S$ such that $x(t) \in S$ for all $t$, which we assume from now on. Let $\bar{x} = \bar{x}_S$. Notice that (65) again holds, for all $x \in \mathbb{R}^n_+$. Let us now specialize to systems (5). We claim that

$$\langle \vec{\rho}(x) - \vec{\rho}(\bar{x}), f(x) \rangle < 0 \tag{69}$$

(i.e., property (62) for (5)) is valid for all $x \in S \cap \mathbb{R}^n_+$, $x \neq \bar{x}$. Since $c(\bar{x}) > 0$, by (65) all we need to show is that $\delta(x, \bar{x}) \neq 0$. But the only way that $\delta(x, \bar{x})$ could vanish is if $x$ is an equilibrium (cf. Corollary 5.3), and uniqueness of equilibria in $S \cap \mathbb{R}^n_+$ then gives $x = \bar{x}$. So, Part 1 in Lemma 8.4 gives $x(t) \to E$. ∎

### 8.2.3 Proof of Theorem 4

Let $S$ be a positive class, and pick $\bar{x} = \bar{x}_S$. Define:

$$\delta_S(x) := \frac{1}{4} c(\bar{x})^2\, \delta(x, \bar{x})$$

restricted to $x \in S \cap \mathbb{R}^n_+$. Since $c(\bar{x}) > 0$ and $\delta(x, \bar{x}) = 0$ if and only if $x$ is an equilibrium, which by uniqueness means $x = \bar{x}$, we have $\delta_S(x) > 0$ for all $x \neq \bar{x}$. take now any collection of functions $\{\Delta_{ij}\}$ such that (22) holds.

Notice that (65) again holds, for all $x \in \mathbb{R}^n_+$. On the other hand, for $x \in S \bigcap \mathbb{R}^n_+$:

$$\begin{aligned}
\langle \vec{\rho}(x) - \vec{\rho}(\bar{x}), g(x) \rangle &= \sum_{i=1}^m \sum_{j=1}^m \Delta_{ij}(x) \langle \vec{\rho}(x) - \vec{\rho}(\bar{x}), b_i - b_j \rangle \\
&\leq \sqrt{\delta_S(x)}\sqrt{\delta(x,\bar{x})} \\
&\leq \frac{1}{2} c(\bar{x}) \delta(x,\bar{x}),
\end{aligned}$$

where the first inequality follows from the Cauchy-Schwartz inequality, the estimate (22), and the definition of $\delta$. Therefore $\langle \vec{\rho}(x) - \vec{\rho}(\bar{x}), f^*(x) \rangle$ equals:

$$\langle \vec{\rho}(x) - \vec{\rho}(\bar{x}), f(x) \rangle + \langle \vec{\rho}(x) - \vec{\rho}(\bar{x}), g(x) \rangle \leq -\frac{1}{2} c(\bar{x}) \delta(x,\bar{x}) \tag{70}$$

for all $x \in S \bigcap \mathbb{R}^n_+$, and this expression is negative when $x \neq \bar{x}$. Thus, hypothesis (62) in Lemma 8.4 holds, with $S$ as given. Part 3 (asymptotic stability of $\bar{x}$) follows from the Lemma.

We pick any $\xi \in S$, and consider the ensuing maximal solution. By Part 1 of Lemma 8.4, the solution is defined for all $t \geq 0$ (proving Part 2 in Theorem 4) and $x(t) \to E_0 \bigcup \{\bar{x}\}$ as $t \to +\infty$. Since $E_0$ and $\{\bar{x}\}$ are at positive distance, this means that either $x(t) \to E_0$ or $x(t) \to \bar{x}$. In the first case, $S$ being closed implies that $x(t) \to E_0 \bigcap S$. Thus, if $E_0 = E_0^*$ does not intersect $S$, the only possibility is that $x(t) \to \bar{x}$. Conversely, if $S \bigcap E_0 \neq \emptyset$ then $\bar{x}_S$ is not globally asymptotically stable relative to $S$. This is clear, since if $\zeta \in E_0 \bigcap S \neq \emptyset$ then $\zeta$ (being an equilibrium) is not attracted to $\bar{x}$. This proves Part 4 of Theorem 4. ∎

Observe $\delta(x, \bar{x}) = \delta(x, \pi(x))$ is of class $\mathcal{C}^k$ if hypothesis ($H_k$) holds. One may pick easily a smooth lower bound for the function $c$, and in that manner obtain a class $\mathcal{C}^k$ function $\delta_S$.

### 8.2.4 Proof of Parts b and c in Theorem 2 and Lemma 2.2

Pick any positive class $S$. We let $\delta_S$ be any function as in the statement of Theorem 4; as $\Delta_{ij}(x) \equiv 0$, the hypotheses of that theorem are verified. In this manner, all conclusions in Theorem 2 as well as Lemma 2.2 are established. ∎

**Remark 8.5** A different approach to the proof of Theorem 1, not using LaSalle invariance, can be based upon the fact that $V(x(t))$ can be shown to decrease strictly along *every* trajectory (not merely in the positive orthant), and this proof applies as well in the multiple-linkage class case described in Remark 5.5. (The fact that all trajectories approach equilibria, in the multiple-linkage class case, may also be proved directly. The LaSalle argument proceeds in basically the same manner; the critical step is to study those trajectories that remain in a level set of the Lyapunov function and start on the boundary. In the case of a single class, such trajectories are either equilibria or enter immediately the interior $\mathbb{R}^n_+$ (Corollary 7.5). In the general case, there may be a block which is "turned off", that is, $f_i(x(t)) \equiv 0$ for some $i \in \{1, \ldots, L\}$. In that event, we are reduced to studying a system with $L-1$ classes, and the same Lyapunov function (restricted to a possibly smaller state space, after dropping those state variables which only appear in the "off" reaction) can be used inductively. □

## 8.3 A Remark on Exponential Stability

Suppose that hypothesis ($H_1$) holds. Recalling the definition (63) of the Lyapunov function $V$, and the form of its gradient given by (59), $(\vec{\rho}(x) - \vec{\rho}(\bar{x}))'$, we know that the Hessian of $V$ at $\bar{x}$ is given by $Q = \mathrm{diag}(\rho'_1(\bar{x}_1), \ldots, \rho'_n(\bar{x}_n))$ and is therefore nonsingular. Since the gradient of $V$ vanishes at $x = \bar{x}$, we have a Taylor expansion $V(x) = (x-\bar{x})' Q (x-\bar{x}) + o(|x-\bar{x}|^2)$, and thus, for all $x$ in some neighborhood of $\bar{x}$, $c_1 |x-\bar{x}|^2 \leq V(x) \leq c_2 |x-\bar{x}|^2$ for some positive constants $c_1$ and $c_2$. If we let $d_2 = d_2(\bar{x})/2$ in (36), we also have that $\delta(x,\bar{x}) \geq d_2 |x - \bar{x}|^2$ for all $x$ in a neighborhood of $\bar{x}$. So inequality (70) gives us that, along trajectories in the class containing $\bar{x}$, as long as $x(t)$ is near $\bar{x}$ we have

$$\frac{dV(x(t))}{dt} = \langle \vec{\rho}(x) - \vec{\rho}(\bar{x}), f^*(x) \rangle \leq -c\,V(x(t)),$$

where $c = c(\bar{x})d_2 c_1/2$. Integrating, $V(x(t)) \leq e^{-ct}V(x(0))$, so we obtain an estimate $|x(t) - \bar{x}|^2 \leq (c_2/c_1)e^{-ct}|x(0) - \bar{x}|^2$ for all trajectories of system (21) which start near $\bar{x}$. In other words, (relative) stability is in fact *exponential*. ∎

# References

[1] G. Bastin, "Issues in modelling and control of mass balance systems," in *Stability and Stabilization of Nonlinear Systems* (D. Aeyels, F. Lamnabhi-Lagarrigue, and A.J. van der Schaft, eds), *Springer-Verlag*, Berlin, 1999, pp. 53-74.

[2] G. Bastin and D. Dochain, *On-line Estimation and Adaptive Control of Bioreactors*, Elsevier, Amsterdam, 1990.

[3] G. Bastin and J. Levine, "On state accessibility in reaction systems," *IEEE Trans. Autom. Control* **38**(1992), pp. 733-742.

[4] G. Bastin and J.F. van Impe, "Nonlinear and adaptive control in biotechnology: A tutorial," *European J. Control* **1**(1995), pp. 1-37.

[5] D.S. Bernstein and S.P. Bhat, "Nonnegativity, Reducibility, and Semistability of Mass Action Kinetics," *Proc. Conf. Dec. Contr.*, pp. 2206-2211, Phoenix, AZ, December 1999.

[6] M. Chaves and E.D. Sontag, "State-estimators for chemical reaction networks of Feinberg-Horn-Jackson zero deficiency type," submitted.

[7] I.R. Epstein and J.A. Pojman, *An Introduction to Nonlinear Chemical Dynamics. Oscillations, Waves, Patterns, and Chaos*, Oxford University Press, New York, 1998.

[8] J. Fabrice, *Dynamics and Robust Nonlinear Control of Stirred Tank Reactors*, Thése de Docteur en Sciences Appliquées, Université catholique de Louvain, Belgium, 1996.

[9] M. Feinberg, "Chemical reaction network structure and the stabiliy of complex isothermal reactors - I. The deficiency zero and deficiency one theorems," Review Article 25, *Chemical Engr. Sci.* **42**(1987), pp. 2229-2268.

[10] M. Feinberg, "The existence and uniqueness of steady states for a class of chemical reaction networks," *Archive for Rational Mechanics and Analysis* **132**(1995), pp. 311-370.

[11] M. Feinberg, "Lectures on Chemical Reaction Networks," 4.5 out of 9 lectures delivered at the Mathematics Research Center, University of Wisconsin, Fall, 1979.

[12] M. Feinberg, "Mathematical aspects of mass action kinetics," in *Chemical Reactor Theory: A Review* (L. Lapidus and N. Amundson, eds.), Prentice-Hall, Englewood Cliffs, 1977.

[13] O. Françoisse, *Modélisation et Commande Adaptative de Réacteurs Tubulaires*, Thése de Docteur en Sciences Appliquées, Université catholique de Louvain, Belgium, 1993.

[14] G.R. Gavalas, *Nonlinear Differential Equations of Chemically Reacting Systems*, Springer-Verlag, Berlin, 1968.

[15] J.K. Hale, *Ordinary Differential Equations*, Wiley, New York, 1980.

[16] F.J.M. Horn and R. Jackson, "General mass action kinetics," *Arch. Rational Mech. Anal.* **49**(1972), pp. 81-116.

[17] F.J.M. Horn, "The dynamics of open reaction systems," in *Mathematical aspects of chemical and biochemical problems and quantum chemistry (Proc. SIAM-AMS Sympos. Appl. Math., New York, 1974)*, pp. 125-137. SIAM-AMS Proceedings, Vol. VIII, Amer. Math. Soc., Providence, 1974.

[18] J.A. Jacquez and C.P. Simon, "Quanlitative theory of compartmental systems," *SIAM Review* **35**(1993), pp. 43-79.

[19] T.W. McKeithan, "Kinetic proofreading in T-cell receptor signal transduction," *Proc. Natl. Acad. Sci. USA* **92**(1995), pp. 5042-5046.

[20] D. Siegel and S.F. Chen, "Global stability of deficiency zero chemical networks," *Canadian Applied Math Quarterly* **2**(1994), pp. 413-434.

[21] E.D. Sontag, *Mathematical Control Theory: Deterministic Finite Dimensional Systems, Second Edition*, Springer-Verlag, New York, 1998.

[22] E.D. Sontag, "Global stability of McKeithan's kinetic proofreading model for T-cell receptor signal transduction," *Mathematics ArXiv*, paper math.DS/9912237, 1999.

[23] A. Vol'pert and S. Hudjaev, *Analysis in Classes of Discontinuous Functions and Equations of Mathematical Physics*, Marinus Nijhoff, Dordrecht, 1985.

[24] C. Wofsy, "Modeling receptor aggregation," series of lectures at the *Workshop on Mathematical Cellular Biology*, Pacific Institute for the Mathematical Sciences, UBC, Vancouver, August 1999.